\documentstyle[11pt,french,amssymb,amsfonts,amsmath,xypic]{article}
\def\tanmoytranslate#1{\catcode`#1=\active\mytanmoytranslate}
\def\mytanmoytranslate#1#2{\def#1{#2}}
\tanmoytranslate    \ä ä{\"a}
\tanmoytranslate    \à à{\`a}
\tanmoytranslate    \á á{\'a}
\tanmoytranslate    \ã ã{\~a}
\tanmoytranslate    \â â{\^a}
\tanmoytranslate    \ë ë{\"e}
\tanmoytranslate    \è è{\`e}
\tanmoytranslate    \é é{\'e}
\tanmoytranslate    \ê ê{\^e}
\tanmoytranslate    \ï ï{\"\i}
\tanmoytranslate    \ì ì{\`\i}
\tanmoytranslate    \í í{\'\i}
\tanmoytranslate    \î î{\^\i}
\tanmoytranslate    \ö ö{\"o}
\tanmoytranslate    \ò ò{\`o}
\tanmoytranslate    \ó ó{\'o}
\tanmoytranslate    \õ õ{\~o}
\tanmoytranslate    \ô ô{\^o}
\tanmoytranslate    \ü ü{\"u}
\tanmoytranslate    \ù ù{\`u}
\tanmoytranslate    \ú ú{\'u}
\tanmoytranslate    \û û{\^u}
\tanmoytranslate    \Ä Ä{\"A}
\tanmoytranslate    \À À{\`A}
\tanmoytranslate    \Á Á{\'A}
\tanmoytranslate    \Ã Ã{\~A}
\tanmoytranslate    \Â Â{\^A}
\tanmoytranslate    \Ë Ë{\"E}
\tanmoytranslate    \È È{\`E}
\tanmoytranslate    \É É{\'E}
\tanmoytranslate    \Ê Ê{\^E}
\tanmoytranslate    \Ï Ï{\"I}
\tanmoytranslate    \Ì Ì{\`I}
\tanmoytranslate    \Í Í{\'I}
\tanmoytranslate    \Î Î{\^I}
\tanmoytranslate    \Ö Ö{\"O}
\tanmoytranslate    \Ò Ò{\`O}
\tanmoytranslate    \Ó Ó{\'O}
\tanmoytranslate    \Õ Õ{\~O}
\tanmoytranslate    \Ô Ô{\^O}
\tanmoytranslate    \Ü Ü{\"U}
\tanmoytranslate    \Ù Ù{\`U}
\tanmoytranslate    \Ú Ú{\'U}
\tanmoytranslate    \Û Û{\^U}
\tanmoytranslate    \ñ ñ{\~n}
\tanmoytranslate    \Ñ Ñ{\~N}
\tanmoytranslate    \ç ç{\c c}
\tanmoytranslate    \Ç Ç{\c C}
\tanmoytranslate    \ß ß{\ss}
\tanmoytranslate    \Æ Æ{\AE}
\tanmoytranslate    \æ æ{\ae}
\tanmoytranslate    \Å Å{\AA}
\tanmoytranslate    \å å{\aa}
\tanmoytranslate    \© ©{\copyright}
\tanmoytranslate    \£ £{\pounds}
\tanmoytranslate    \¶ ¶{\P}
\tanmoytranslate    \§ §{\S}
\tanmoytranslate    \¿ ¿{?`}
\tanmoytranslate    \¡ ¡{!`}

\begin{document}

\author{Gentiana Danila}
\title{Sur la cohomologie d'un fibré tautologique sur le schéma de
  Hilbert d'une surface  }
\date{1 Avril 1999}
\maketitle
\newtheorem{theor}{Théorème}[section]
\newtheorem{prop}[theor]{Proposition}
\newtheorem{cor}[theor]{Corollaire}
\newtheorem{lemme}[theor]{Lemme}
\newtheorem{slemme}[theor]{Sous-lemme}
\newtheorem{defi}[theor]{Définition}
\newtheorem{defiprop}[theor]{Définition-Proposition}
\newtheorem{conj}[theor]{Conjecture}
\newtheorem{rem}[theor]{Remarque}
\newtheorem{rems}[theor]{Remarques}
\newtheorem{nota}[theor]{Notations}
\newtheorem{rappel}[theor]{Rappel}

\def\fs{{faisceau }}  
\def\fx{{faisceaux }}
\def\alg{{algébrique }}
\def\algs{{algébriques }}
\def\th{{théorème }}
\def\rep{{représentation }}
\def\reps{{représentations }}
\def\irr{{irréductible }}
\def\irrs{{irréductibles }}
\def\coh{{cohomologie }}
\def\co{{cohérent }}
\def\fib{{fibré }}
\def\fibs{{fibrés }}
\def\mor{{morphisme }}
\def\isom{{isomorphisme }}
\def\iso{{isomorphisme }}
\def\mors{{morphismes }}
\def\sur{{surjectif }}
\def\diff{{différentielle }}
\def\diffs{{différentielles }}
\def\inve{{inversible }}
\def\sct{{section }}
\def\hol{{holomorphe }}
\def\hols{{holomorphes }}
\def\app{{application }}
\def\apps{{applications }}

\def\proj{{\mathbb{P}}}
\def\hilx{X^{\mbox{}^{[m]}}}
\def\hilxp{X^{\mbox{}^{[m+1]}}}
\def\hilxmm{X^{\mbox{}^{[m,m+1]}}}
\def\hilxs{X_*^{\mbox{}^{[m]}}}
\def\hilxx{\hilx\times X}
\def\xm{X^m}
\def\xms{X^m_*}
\def\xmp{X^{m+1}}
\def\smx{\es^m(X)}
\def\smxp{\es^{m+1}(X)}
\def\smxs{\es^m_*(X)}
\def\blup{{\rm{Bl}}_{\Xi_m}(\hilx\times X)}
\def\blzx{{\rm{Bl}}_{Z}(X)}
\def\blym{{\rm{Bl}}_{Y}(M)}
\def\blim{{\rm{Bl}}_{I}(M)}
\def\pim{\proj_M(I)}
\def\pxiz{\proj_X(\I_Z)}
\def\piy{\proj(I_Y)}
\def\pxi{\proj(\ix)}
\def\zu{Z_{univ}}
\def\zup{Z^{\prime}_{univ}}
\def\zupp{Z^{\prime\prime}_{univ}}
\def\xiu{\xi_{univ}}
\def\xipu{\xi^{\prime}_{univ}}
\def\xippu{\xi^{\prime\prime}_{univ}}
\def\etau{\eta_{univ}}
\def\etapu{\eta^{\prime}_{univ}}
\def\etappu{\eta^{\prime\prime}_{univ}}
\def\zs{Z}
\def\xis{\xi}
\def\xips{\xi^{\prime}}
\def\etas{\eta}
\def\xim{\Xi_{m}}
\def\ximp{\Xi_{m+1}}
\def\pp{\proj_2}
\def\ppm{\pp^{\mbox{}^{[m]}}}

\def\de{{\mathfrak{d}}^A_m}
\def\deo{{\mathfrak{d}}^{\ox}_m}
\def\dep{{\mathfrak{d}}^A_{m+1}}
\def\dem{{\mathfrak{d}}^A_{m,1}}
\def\lm{L^{\mbox{}^{[m]}}}
\def\om{\O^{\mbox{}^{[m]}}}
\def\lam{(L\tens A)^{\mbox{}^{[m]}}}
\def\lmp{L^{\mbox{}^{[m+1]}}}
\def\O{{\mathcal O}}
\def\ox{\omega_X}
\def\oy{\omega_Y}
\def\oxy{\omega_{X/Y}}
\def\iu{I_{univ}}
\def\ix{{{I}}_{\xim}}
\def\oxm{\O_{\xim}}
\def\pix{\proj(\ix)}
\def\pa{\proj(\A)}
\def\pf{\proj(\F)}
\def\dmp{{\mathcal{D}}_{m+1}}
\def\dm{{\mathcal{D}}_{m}}

\def\sigm{{\mathfrak{S}}_m}
\def\comp{{\Bbb C}}
\def\sigmp{{\mathfrak{S}}_{m+1}}

\def\rightto#1{\smash{\mathop{\longrightarrow}\limits^{#1}}}
\def\downto#1{\Big\downarrow\rlap{$\vcenter{\hbox{$\scriptstyle#1$}}$}}
\def\surto{\twoheadrightarrow}

\def\tens{\otimes}

\def\L{{\mathcal L}}
\def\B{{\mathcal B}}
\def\A{{\mathcal A}}
\def\K{{\mathcal K}}
\def\has{{\mathcal H}}
\def\I{{\mathcal I}}
\def\F{{\mathcal F}}
\def\N{{\mathcal N}}

\def\uHom{\underline{{\rm Hom}}}

\def\etab{\bar{\eta}}
\def\tL{\widetilde{L}}
\def\tV{\widetilde{V}}
\def\tI{\widetilde{I}}
\def\tW{\widetilde{W}}
\def\tK{\widetilde{K}}
\def\tF{\widetilde{\F}}
\def\tM{\widetilde{M}}
\def\tD{\widetilde{D}}
\def\tmL{\widetilde{\L}}

\def\M{{\bf M}}

\def\es{{\rm S}}
\def\H{{\rm H}}
\def\Stab{{\rm Stab\,} }
\def\Ker{{\rm Ker}\, }
\def\coker{{\rm coker}\, }
\def\dim{{\rm dim}\, }
\def\im{{\rm im}\, }
\def\codim{{\rm codim}\, }
\def\supp{{\rm supp}\, } 
\def\Hom{{\rm Hom}\, }
\def\Pic{{\rm Pic}\, }
\def\Tor{{\rm Tor}\, }

{\bf Abstract~:} {\small We compute the cohomology spaces for the
tautological bundle tensor the determinant bundle on the punctual
Hilbert scheme $\hilx$ of a smooth projective surface $X$ on
$\comp$. We show that for $L$ and $A$ invertible vector bundles on
$X$, and $\omega_X$ the canonical bundle of $X$, if $\ox^{-1}\otimes
A$, $\ox^{-1}\otimes L$  and $A$ are ample
vector bundles, then the higher cohomology spaces on
$\hilx$ of the tautological bundle associated to $L$ tensor the
determinant bundle associated to $A$ vanish, and the space of global
sections is computed in terms of $\H^0(A)$ and $\H^0(X, L\otimes
A)$. This result is motivated by the computation of the space of
global sections of the determinant bundle on the moduli space of rank
$2$ semi-stable sheaves on the projective plane, supporting Le
Potier's Strange duality conjecture on the projective plane.} 



Running heads: Fibré tautologique sur le schéma de Hilbert d'une surface

\section{Introduction}
\label{sec:int}

Soit $X$ une surface projective lisse sur $\comp$, $L$ et $A$ deux
fibrés inversibles sur $X$, $\omega_X$ le fibré canonique de $X$. Pour
tout entier $m$, on note $\hilx$ le schéma de Hilbert qui paramètre les
sous-schémas finis de $X$ de longueur $m$. Il est lisse et projectif
de dimension $2m$. On utilise le diagramme
$$\diagram
&\xm\dto^{q}\\
\hilx\rto^{HC}&\smx
\enddiagram
$$
pour construire sur $\hilx$ un \fib $\de$ associé à $A$. Le morphisme
$HC$ est le morphisme de Hilbert-Chow, qui associe à un sous-schéma
$Z\subset X$ le
cycle $\sum_{x\in X} lg(Z_x)x$, où $Z_x$ est la composante de $Z$
passant par $x$, et $lg(Z_x)$ la longueur de $Z_x$. Le \mor  $q$ est
le quotient par l'action naturelle du groupe symétrique
$\sigm$. On définit $\de=HC^*((A\boxtimes A\boxtimes\cdots\boxtimes
A)^{\sigm})$. Via
l'identification de $\hilx$ à l'espace de modules des \fx sans torsion
de rang $1$ sur $X$, on retrouve le \fib déterminant sur $\hilx$
associé à $A$. Soit
$\xim\subset\hilx\times X$, le schéma universel des couples $(Z,x)$ tels
que $x\in Z$. Il est fini de degré $m$ et plat au-dessus de
$\hilx$. On note $p_1$ et $p_2$ les  deux projections~:
$$\diagram
\hilxx\rto^{\ \  p_2}\dto_{p_1}& X\\
\hilx
\enddiagram
$$
aussi bien que leurs restrictions à $\Xi_m$, et $\lm$ le \fib $p_{1*}(\O_{\Xi_m}\otimes p_2^*(L))$.

Le but de cet article est de démontrer le
\begin{theor}
\label{the:1}
Si le \fib $\ox^{-1}\otimes A^{\tens k}$ est ample pour $1\le k\le m$ et
$\ox^{-1}\otimes L\otimes A^{\tens k}$ est ample pour $1\le k\le m$, alors

i) $\H^q(\hilx, \lm\otimes\de)=0$ pour $q>0$;

ii) $\H^0(\hilx, \lm\otimes\de)\simeq\es^{m-1}(\H^0(A))\otimes\H^0(X, L\otimes A).$
\end{theor}

Ce \th admet plusieurs corollaires, dont un en particulier est utile
dans l'article \cite{D}, pour le calcul de l'espace des sections du
\fib déterminant sur l'espace de modules des faisceaux semi-stables de
rang $2$ sur le plan projectif. Ce calcul fournit des exemples pour la
conjecture de Le Potier de ``Dualité étrange'' sur  le plan
projectif. On s'appuyera sur les travaux
\cite{E-L}, \cite{E-S}, \cite{Lehn}, \cite{Tik}, \cite{Che} de
G. Ellingsrud, S. Str\o mme, M. Lehn, A. Tikhomirov, J. Cheah.

\begin{cor}
  Les assertions i) et ii) sont vraies si les fibrés $\ox^{-1}\otimes
  A$, $\ox^{-1}\otimes L$ et $A$ sont amples.
\end{cor}

On utilise le fait que le produit tensoriel de deux fibrés amples est
ample.

\begin{cor}
  Si les fibrés $\ox^{-1}$ et $\ox^{-1}\otimes
  L$ sont amples alors

i) $\H^q(\hilx,\lm)=0$ pour $q>0$;

ii) $\H^0(\hilx,\lm)\simeq\H^0(X,L)$.
\end{cor}

On applique le \th pour $A=\O_X$.

En particulier si $X$ est le plan projectif complexe $\pp$, $L=\O_{\pp}(3)$ et
  $A=\O_{\pp}(1)$, ou si $L=\O_{\pp}(-1)$ et
  $A=\O_{\pp}(3)$ alors

i) $\H^q(\ppm,\lm\tens\de)=0$ pour $q>0$;

ii) $\H^0(\ppm,\lm\tens\de)\simeq\es^{m-1}(\H^0(A))\tens\H^0(\pp,L\tens A)$.

Pour $L=\O_{\pp}(3)$, $A=\O_{\pp}(1)$, c'est le résultat
utilisé dans \cite{D}.

Les cas $q=0$ et $q=1$ se traitent directement et leur démonstration est
donnée dans les préliminaires. Pour $q>1$ on utilise une récurrence
sur $m$. Il nous faut introduire la variété auxiliaire $\hilxmm$
d'incidence. C'est le sous-schéma fermé de
$\hilxp\times\hilx$ donné par $\hilxmm=\{(Z,\xi)| \xi\subset Z\}$. Il
est connu que $\hilxmm$ est lisse et irréductible de dimension $2m+2$
(\cite{Che}, \cite{Tik}). Il y a des \mors évidents $p_m:\hilxmm\to\hilx$ et
$p_{m+1}:\hilxmm\to\hilxp$ induits par les projections. Il y a aussi
un \mor naturel $q:\hilxmm\to X$ qui envoie une paire $(Z,\xi)$ sur
l'unique point $\eta$ où les schémas $\xi$ et $Z$ diffèrent (schématiquement). Dans la
démonstration du \th \ref{the:1} on utilisera le passage de
$\hilxx$ à $\hilxp$ à travers $\hilxmm$ et les \mors $p_{m+1}$
et $\phi=(p_m,q)$:
$$\diagram
\hilxmm\rto^{\phi}\dto_{p_{m+1}}&\hilx\times X\\
\hilxp&.
\enddiagram
$$

\section{Préliminaires}
\label{sec:pre}

Notations~: Le corps de base est $\comp$. Pour un espace vectoriel $V$ nous noterons $\proj(V)$  l'espace projectif de
Grothendieck des espaces vectoriels
quotients de dimension $1$. Par variété algébrique on entend schéma de
type fini sur $\comp$, séparé; les points considérés sont toujours les
points fermés.

\subsection{Cohomologie à support}
\label{sec:coh}

Pour un \fs abélien $F$ sur un espace topologique $X$, les espaces de \coh à support et les \fx  de cohomologie à support
dans un  fermé $Y$ de $X$: $\H^i_Y(F)$,   et $\has^i_Y(F)$, sont  définis
dans le premier paragraphe de  \cite{grot}, et ils sont reliés par une suite spectrale 
$$E^{p,q}_2=\H^p(X,\has^q_Y(F))\Rightarrow\H^{n}_Y(F)$$ 
d'aboutissement $\H^n_Y(F)$ en degré $n=p+q$ (\cite{grot}, prop 1.4, p.5). 

On associe à $Y\subset X$ une suite exacte de \coh locale
(\cite{grot}, cor. 1.9, p.9):
$$\cdots\to \H^i_Y(X,F)\to \H^i(X,F)\to \H^i(X\setminus Y,F)\to
\H^{i+1}_Y(X,F)\to\cdots $$

Nous utiliserons le
lemme suivant~:

\begin{lemme}(\cite{grot}, thm. 3.8, p.44):
\label{lem:2.1}
Soient $X$ une variété
\alg  lisse, $F$ un \fs de $\O_X$-modules localement libre sur $X$  et $Y\subset X$ un  fermé de $X$. Pour un
entier $n$ donné les assertions suivantes sont équivalentes: 
\begin{itemize}
\item{\rm (i)} $\codim Y\ge n$ 
\item{\rm (ii)} pour tout $ i<n$, $\has^i_Y(F)=0$. 
\end{itemize}
\end{lemme}

On  déduit  de la suite spectrale ci-dessus que si $\codim Y\ge n$ 
\begin{equation}
\H^i_Y(F)=0 {\text{ pour tout }} i<n  {\text{ et }}
\H^n_Y(F)=\H^0(\has^n_Y(F)).  
\label{equ:2}
\end{equation}

\subsection{Calculs d'invariants}
\label{sec:cal}

On considère un ensemble fini $I$ muni d'une action transitive d'un
groupe fini $G$. Soit $Y$ une variété sur laquelle $G$ agit à
gauche. Considérons  pour chaque $i\in I$ un \fib $L_i$ sur $Y$ de fa\c con 
qu'on ait un \iso canonique $h_g:g^*(L_i)\simeq 
L_{g^{-1}(i)}$ pour tout $i\in I$ et $g\in G$, et  pour tous $g,g'\in G$, $h_g\circ h_{g'}=h_{gg'}$. (En 
particulier pour tout $g$ dans $\Stab\{i\}$, le stabilisateur de $i$, on a  $g^*(L_i)\simeq 
L_i$).  On a alors un diagramme commutatif:
$$\diagram
L_{g^{-1}(i)}\dto\rto^{h_{g}}&L_i\dto\\
Y\rto^{g}&Y\enddiagram$$

On considère l'espace vectoriel des sections $M_i=\H^0(L_i)$, et la 
somme directe $M=\oplus_{i\in I}M_i$ (espace vectoriel des familles 
$s=(s_i)_i$ avec $s_i\in M_i$).

L'\iso $h_{g}$ induit un \iso $\lambda_{g}:M_i\to 
M_{g(i)}$ en posant pour $x\in Y$
$$\lambda_{g}(s)(x)=h_{g}s(g^{-1}(x)).$$

On peut facilement vérifier que $\lambda_{gg'}=\lambda_{g}\lambda_{g'}$. En 
particulier, ceci fournit une 
action à gauche du stabilisateur de $i$ sur $M_i$. On définit aussi 
une action à gauche de $G$ sur $M$ en posant 
$g(s)_i=\lambda_{g}(s_{g^{-1}(i)})$. Le lemme 
suivant  est l'ingrédient essentiel pour
les calculs d'invariants sur $\xm$:

\begin{lemme}
\label{lemm:2.1}
Soit $i\in I$. La projection 
$pr_i:M\to M_i$ induit un \iso $M^G\to M_i^{\Stab\{i\}}$.
\end{lemme}

\par{\bf Preuve du lemme~:}
 D'après la définition de l'action  l'image par 
$pr_i$ des invariants de $M$ par $G$ est contenue dans le 
sous-espace des invariants de $M_i$ par $\Stab\{i\}$.

Montrons l'injectivité: si $u\in M^{G}$, pour tout élément
$g$ tel que $g(i)=j$, on a $u_j=\lambda_{g}(u_i)$, 
ce qui montre que $u$ est déterminé par $u_i$.

Montrons la surjectivité: soit $v\in M_i$, invariant par $\Stab\{i\}$. 
On définit $u\in M$ par la formule $u_j=\lambda_{g}(v)$ où 
$g$ est un élément tel que $g(i)=j$. Il faut 
bien sûr vérifier que ceci ne dépend pas du choix de $g$; mais ceci résulte de l'hypothèse que $v$ est invariant. 
En effet, si $g'$ est une autre élément de $G$ tel que 
$g'(i)=j$, l'élément $g^{-1}g'$ appartient à 
$\Stab\{i\}$ et on a 
$\lambda_{g'}(v)=\lambda_{g}\lambda_{g^{-1}g'}(v)=
\lambda_{g}(v)$. Il est clair que l'élément $u$ construit est 
invariant par $G$.$\Box$

\subsection{Les espaces $\H^0$ et $\H^1$}
\label{sec:les}

Soient, comme dans l'introduction, $L$ et $A$ deux \fibs inversibles sur
une surface projective lisse $X$.
\begin{prop}
\label{pro:1}
Si $\ox^{-1}\otimes A$  et
$\ox^{-1}\otimes L\otimes A$ sont  amples alors

i) $\H^1(\hilx, \lm\otimes\de)=0;$

ii) $\H^0(\hilx, \lm\otimes\de)\simeq\es^{m-1}(\H^0(A))\otimes\H^0(X, L\otimes A).$
\end{prop}

La méthode utilisée consiste à calculer ces espaces sur un grand
ouvert $\hilxs$ de $\hilx$ et à utiliser des résultats de \coh à
support.

L'ouvert $\hilxs$ est formé par les schémas avec au 
plus un point multiple, qui soit double, soit les schémas dont le 
cycle correspondant est $x_1+x_2+\cdots+x_m$ ou $2x_1+x_3+x_4+\cdots+x_m$ avec $x_i$ distincts.
On note $\smxs$ l'ouvert des cycles de cette forme. Rappelons que  $q:\xm\to\smx$ est le 
quotient  de $\xm$ sous l'action du groupe symétrique $\sigm$.
 L'avantage d'utiliser $\hilxs$ 
est qu'on peut le décrire comme quotient $p$ 
de l'éclaté $B$ de $\xms=q^{-1}(\smxs)$ selon la réunion $D$ 
des diagonales $\Delta_{ij}=\{(x_1,\cdots,x_m)\in\xms|x_i=x_j\}$ pour 
$i<j$, disjointes dans $\xms$. On note $\rho$ cet éclatement. On a un 
diagramme commutatif:
$$\diagram
B  \rto^{\rho}\dto^{p} & \xms\dto^{q} \\
\hilxs \rto^{\pi} &\smxs
\enddiagram$$

On va montrer comment, à l'aide de cette description, on peut ramener 
les calculs de la \coh des fibrés sur $\hilxs$ à des calculs des 
invariants de la \coh de certains \fx sur $\xms$.

\par{\bf Preuve de la proposition:}
 Sur $B$, le
diviseur exceptionnel $E$ se décompose en  
composantes disjointes $E=\bigcup_{i<j}E_{i,j}$. Alors le schéma 
universel $\Xi_B\subset B\times X$, paramétré par $B$, a $m$ 
composantes irréductibles $\Xi_i$ et la projection 
$p_1:\Xi_i\bigcap\Xi_j\to E_{i,j}$ est un isomorphisme. On en déduit 
une suite exacte sur $B\times X$:
\begin{equation}
0\to\O_{\Xi_B}\to\oplus_i\O_{\Xi_i}\to\oplus_{i<j}\O_{E_{i,j}}\to 0 
\label{s6}
\end{equation}
Comme, par changement de base, $p^{*}(\lm)=p_{1*}(\O_{\Xi_B}\tens 
p_2^*(L))$, on a, après tensorisation par 
$p^*_2(L)$ de la suite (\ref{s6}) et image directe par $p_1$
(qui est un \mor fini en restriction au schéma universel), une suite sur $B$:
\begin{equation}
\label{s6bis}
0\to
p^{*}(\lm)\to\oplus_ip_i^*(L)\stackrel{\epsilon}{\to}\oplus_{i<j}p^*_{i,j}(L_\Delta)\to
0.
\end{equation}
La projection $p_i$ désigne aussi bien la $i$-ème projection $\xms\to X$ que 
sa composée avec $\rho:B\to X$; de même pour $p_{i,j}:\xms\to 
X\times X$ et $p_{i,j}:B\to X\times X$. Le sous-schéma $\Xi_i$ est 
l'image réciproque de la diagonale $\Delta$ de $X\times X$ par 
l'application $(p_i,id_X)$. Le \fib $L_\Delta$ est l'image réciproque 
de $L$ par l'une des projections de la diagonale de $X\times X$ sur $X$, 
qui sont des isomorphismes.

Compte-tenu du fait que $p^{*}(\de)=\rho^*(A\boxtimes A\boxtimes\cdots\boxtimes
A)$, il
faut encore tensoriser la suite (\ref{s6bis}) par
$\rho^*(A\boxtimes A\boxtimes\cdots\boxtimes
A)=\otimes^m_{j=1}A_j$ où $A_j=p_j^*(A)$. On introduit aussi les
notations $L_i=p_i^*(L)$ et $L_{i,j}=p^*_{i,j}(L_\Delta)$. Alors on a
la suite exacte
\begin{equation}
\label{s6ter}
0\to
p^{*}(\lm\tens\de)\to\oplus_i(\otimes_{j\ne i}A_j)\tens(A_i\tens L_i) \to\oplus_{i<j}(\otimes_{l\ne i,j}A_j)\tens
p_{i,j}^*((L\tens A)_{\Delta})\to
0.
\end{equation}

Le \fib $\tmL=\oplus_{i<j}(\otimes_{l\ne i,j}A_j)\tens
p_{i,j}^*((L\tens A)_{\Delta})$ a pour support le 
diviseur exceptionnel $E$. Il est l'image réciproque  par $\rho$ du 
\fib $\L=\oplus_{i<j}(\otimes_{l\ne i,j}A_j)\tens
p_{i,j}^*((L\tens A)_{\Delta})$ sur $\xms$, dont le support est 
 $D$, réunion des diagonales $\Delta_{i,j}$. On note $\tW$ et $W$ le \fib $\oplus_i(\otimes_{j\ne i}A_j)\tens(A_i\tens L_i) $ sur $B$ et 
 sur $\xms$ respectivement. Dans ces conditions on remarque que $W$ 
 n'a pas de \coh en degré positif sur $\xm$ (par la formule de
 Künneth et le \th de Kodaira en utilisant l'hypothèse). Puisque le complémentaire de 
 $\xms$ dans $\xm$ est de codimension $4$,  la \coh à support
 dans le complémentaire est nulle en degré $\le 3$ et
 $\H^q(\xm,W)=\H^q(\xms,W)$ pour $q\le 2$ (d'après \ref{sec:coh}).

Le groupe symétrique
$G=\sigm$ opère sur les suites exactes précédentes. La suite (\ref{s6bis}) est
$G$-équivariante. Le \mor $\epsilon$ est donné par 
$(s_i)_i\mapsto(s_i|_{\Delta_{i,j}}-s_j|_{\Delta_{i,j}})_{i,j}$, donc l'action induite sur
$\oplus_{i<j}L_{i,j}$ et sur son image réciproque par $\rho$ est
telle que la transposition $\tau_{i,j}$ change le terme d'indice
$(i,j)$ en son opposé. Aucune section non nulle du \fib
$\oplus_{i<j}L_{i,j}$ n'est donc $G$-équivariante. Considérons une
section invariante $s$ de $\L$. Sa restriction au fermé
$\Delta_{i,j}$, invariant par $\tau_{i,j}$, est
$\tau_{i,j}$-invariante. Comme l'action de
$\tau_{i,j}$ est triviale sur  $\Delta_{i,j}$, $\tau_{i,j}$ agit
trivialement sur $\otimes^m_{j=1}A_j|_{\Delta_{i,j}}$. Par suite
$\tau_{i,j}(s|_{\Delta_{i,j}})=-s|_{\Delta_{i,j}}$, donc
$s|_{\Delta_{i,j}}=0$. Ceci est valable pour tous les $(i,j)$, d'où
$s=0$. Donc  $\L$ n'a pas non
plus de
\coh $G$-équivariante.

Le groupe $G$ 
étant fini, la \coh du \fs des invariants $F^G:=(p_*(F))^G$ sur
$\hilxs$ (ou  $(q_*(F))^G$ sur 
$\smxs$), où $F$ est un $G$-\fs \alg cohérent sur $B$ (ou sur $\xms$)
s'identifie à la  
\coh équivariante de $F$, c'est-à-dire aux invariants de la \coh de
$F$. 

 On écrit alors la suite exacte de \coh
équivariante associée à la suite (\ref{s6ter}):
$$\begin{array}{ccccccccc}
\label{s6q}
0&\to&\H^0(B,p^{*}(\lm\tens\de))^G&\to&\H^0(B,\tW)^G&\to&\H^0(B,\tmL)^G&\to
&\\
&\to&\H^1(B,p^{*}(\lm\tens\de))^G&\to&\H^1(B,\tW)^G&\to&\H^1(B,\tmL)^G&
\to&\cdots
\end{array}$$
Le \mor $\rho$ est le \mor d'éclatement d'une sous-variété lisse d'une
variété lisse. Alors d'après le lemme 3.5 de \cite{SGA-6}, exposé VII,
on obtient 
$$\rho_*(\O_B)=\O_{\xms}$$
et
$$ R^q\rho_*(\O_B)=0 \mbox{ pour }
q>0.$$

On en déduit que $\H^q(B,F)=\H^q(\xms,\rho_*F)$ pour un \fs $F$ sur
$B$ donc $\H^0(B,\tW)^G=\H^0(\xms,W)^G=\H^0(\xm,W)^G$, et
$\H^1(B,\tW)^G=\H^1(\xms,W)^G=\H^1(\xm,W)^G=0.$

L'espace $\H^0(\xm,W)^G=\es^{m-1}(\H^0(A))\tens\H^0(L\tens A)$ se
calcule à l'aide de la section \ref{sec:cal}. Les annulations
$\H^1(B,\tmL)^G=\H^0(B,\tmL)^G=0$ conduisent à
$\H^0(B,p^{*}(\lm\tens\de))^G=\H^0(\hilxs,\lm\tens\de)=\es^{m-1}(\H^0(A))\tens\H^0(L\tens
A)$
et $\H^1(B,p^{*}(\lm\tens\de))^G=\H^1(\hilxs,\lm\tens\de)=0.$
On conclut en tenant compte de la suite exacte de \coh locale associée
à l'ouvert $\hilxs$ dont le complémentaire dans  $\hilx$ est de
codimension $2$.$\Box$

\subsection{Le schéma en espaces projectifs  de Grothendieck associé à un \fs}
\label{sec:surpf}

On rappelle ici des résultats classiques sur le schéma en espaces projectifs de
Grothendieck associé à un faisceau, qu'on peut trouver dans
\cite{EGA-I}, \S 9.7 et \cite{EGA-II}, \S 4. 

Si $M$ est une variété \alg quelconque, et $\F$ est un \fs cohérent sur
$M$, il existe une variété projective $\pf$ munie d'un \mor projectif
$\pi:\pf\to M$, un \fs inversible $\O_{\pf}(1)$ et un \mor \sur
$\pi^*\F\to \O_{\pf}(1)\to 0$, qui résout le problème universel:

Pour tout schéma $S$, tout \mor $f:S\to M$, tout \fib inversible $L$
sur $S$ et tout \mor \sur $f^*\F\to L\to 0$, il existe un \mor unique
$g:S\to\pf$ qui rend commutatif le diagramme
$$\diagram
&\pf\dto^{\pi}\\
S\urto^{g}\rto_{f}&M
\enddiagram$$
et satisfait $L=g^*\O_{\pf}(1)$.

Pour un \fs localement libre de rang $a$, $\F=\A$, $\pa$ est une
fibration avec pour fibre au-dessus de $x\in M$, l'espace projectif
$\proj(\A(x))\simeq\proj^{a-1}$. Lorsque $\F$ n'est pas localement
libre, la solution du problème universel se construit à partir du cas
localement libre de la fa\c con suivante: on considère une résolution
$\B\to\A\to\F\to 0$ de $\F$ par des \fx localement libres $\B$ et
$\A$. Si on note $\pi$ la fibration $\pa\to M$, on obtient la suite
$\pi^*\B\to\pi^*\A\to\pi^*\F\to 0$ sur $\pa$. Le \mor
$\pi^*\A\to\O_{\pa}(1)$ se factorise à travers $\pi^*\A\to\pi^*\F$ si
et seulement si l'application composée 
\begin{equation}
  \label{equ:AC}
  \pi^*\B\to\pi^*\A\to\O_{\pa}(1)
\end{equation}
est nulle. Ainsi $\pf$ est le fermé de $\pa$ donné par l'annulation de
la section $\sigma$ du \fs
$\uHom(\pi^*\B\to\O_{\pa}(1))=\pi^*\B^*\tens_{\O_{\pa}}\O_{\pa}(1)$,
$\O_{\pf}(1)$ est la restriction de $\O_{\pa}(1)$ à ce fermé, et le
\mor $\pi^*\F\to\O_{\pf}(1)$ est obtenu par factorisation de
$\pi^*\A\to\O_{\pa}(1)$.

On aura besoin dans la suite de la relation entre l'éclaté de l'idéal $I$
d'une sous-variété $Y$  de la variété $M$ et le schéma en espaces projectifs associé à cet idéal.
La variété $\blim$ avec le \mor $\rho:\blim\to M$
est objet universel pour la donnée d'une variété $S$ et  d'un \mor
$f:S\to M$ tel que $f^{-1}I\cdot\O_S$ est un \fs inversible. 

\begin{rappel}
\label{pro:soi}
  Soit $M$ une variété, $Y$ une sous-variété de $M$ d'idéal $I$, $\blim$ et $\pim$ les variétés
  associées. Il existe un \mor canonique $b:\blim\to\pim$ au-dessus de $M$.
\end{rappel}

\par{\bf Preuve~:}
Soit $\rho:\blim\to M$ le \mor d'éclatement et $E$ le diviseur
exceptionnel. Par la définition de $\blim$, on a
$\O_{\blim}(-E)=\im(\rho^*I\to\O_{\blim}=\rho^*\O_M)$. Plus exactement
on a sur $M$ la suite exacte~:
$$0\to I\to\O_M\to\O_Y\to 0$$
dont l'image réciproque par $\rho$ sur $\blim$ est
$$\rho^*I\to\O_{\blim}\to\O_E=\rho^*\O_Y\to 0.$$
On dispose également de la suite exacte sur $\blim$ associée au
diviseur $E$~:
$$0\to\O(-E)\to\O_{\blim}\to\O_E\to 0$$
et on a le diagramme commutatif~:
\begin{equation}
  \label{equ:unu}
  {\diagram
&\rho^*I\rto&\O_{\blim}\dto^{\simeq}\rto&\O_E\dto^{\simeq}\rto& 0\\
0\rto&\O(-E)\rto&\O_{\blim}\rto&\O_E\rto& 0.
\enddiagram}
\end{equation}
Par conséquent le \mor $\rho^*I\to\O(-E)$ est \sur, et par la
propriété universelle de $\pim$ on obtient un \mor canonique
$b:\blim\to\pim$ au-dessus de $M$.$\Box$

\subsection{Images directes supérieures}
\label{sec:ima}

Soient $M$ une variété lisse de dimension $d$, $\F$ un \fs \alg
cohérent de rang $r$ et de dimension homologique $\le 1$ (ce qui
revient à l'existence d'une résolution de $\F$ par des \fx localement
libres $\A$ et $\B$: $0\to \B\to \A\to \F\to 0.$) On note $\pi:\pf\to M$ et
on suppose que $\dim\pf=d+r-1$.

\begin{lemme}
\label{lic}
  Dans ces conditions, la variété $\pf$ est localement intersection complète.
\end{lemme}

\par{\bf Preuve du lemme~:}
La variété $\pf$ est un fermé de dimension $d+r-1$ dans la
variété projective lisse $\pa$ (de dimension $d+a-1$, où $a$ est le
rang de $\A$), donnée par
l'annulation d'une section $\sigma$ dans le \fib $\pi^*\B^*(1)$ de
rang $a-r$. On fixe localement sur $\pa$ une base de sections du \fib
$\pi^*\B^*(1)$, $s_1,\ldots, s_{a-r}$. La section $\sigma$ s'écrit $\sigma=\sum
_i f_is_i$ avec $f_i$ des fonctions régulières locales sur
$\pa$. Donc, localement, l'idéal de $\pf$ dans $\pa$ est
$(f_1,\ldots, f_{a-r})$. Comme $\pf$ est de codimension $a-r$ dans
$\pa$, il résulte que $\pf$ est localement intersection
complète. $\Box$

\begin{prop}
\label{dcc}
  Avec les mêmes hypothèses, on obtient

i) $R^q\pi_*\O_{\pf}=\left\{
  \begin{array}{ccc}
0&\mbox{ si }&q>0\\
\O_{M}& \mbox{ si }&q=0
  \end{array}\right.$

ii) $R^q\pi_*\O_{\pf}(1)=\left\{
  \begin{array}{ccc}
0&\mbox{ si }&q>0\\
\F&\mbox{ si }&q=0.
  \end{array}\right.$

\end{prop}
 
\par{\bf Preuve de la proposition~:}

Le lemme \ref{lic} nous permet d'écrire la suite exacte de Koszul sur
$\pa$~:
 \begin{equation}
  \label{equ:(*)}
  \begin{array}{cccccccccc}
0&\to&\Lambda^{a-r}\pi^*\B(r-a)&\to&\Lambda^{a-r+1}\pi^*\B(r-a-1)&\to&\cdots&\to&\Lambda^{2}\pi^*\B(-2)&\to\\
&\to&\pi^*\B(-1)&\to&\O_{\pa}&\to&\O_{\pf}&\to&0.&
  \end{array}
\end{equation}

On utilise le résultat théorique suivant: Si $\pi:Y\to X$ est un \mor
de variétés et 
$$0\to A_n\stackrel{\alpha_n}{\to}
A_{n-1}\stackrel{\alpha_{n-1}}{\to}\cdots\to A_0\to F\to 0$$
est une suite exacte sur $Y$, alors il existe une suite spectrale
$$E_1^{p,q}=R^q\pi_*(A_{-p}) \Longrightarrow  R^{p+q}\pi_*F.$$
La différentielle $\alpha^{p,q}_1:E_1^{p,q}=R^q\pi_*(A_{-p})\to
E_1^{p+1,q}=R^q\pi_*(A_{-p-1})$ est induite par le \mor de \fx
$\alpha_{-p}:A_{-p}\to A_{-p-1}$.
On applique cette suite spectrale pour la résolution de Koszul
(\ref{equ:(*)})~:
$$E_1^{p,q}=R^q\pi_*(\Lambda^{-p}\pi^*\B(p)) \Longrightarrow 
R^{p+q}\pi_*\O_{\pf}.$$
Par convention $\Lambda^{-p}\cdot:=\Lambda^{-p}\cdot$ si $p\le 0$ et
$\Lambda^{-p}\cdot:=0$ si $p>0$. La formule de projection nous donne 
$$R^q\pi_*(\Lambda^{-p}\pi^*\B(p))=\Lambda^{-p}\B\tens_{\O_M}R^q\pi_*(\O_{\pa}(p)).$$
D'après \cite{Hart}, ex. 8.4, p. 253, on a  pour $1-a\le p<0$, et tout
$q$, 
$R^q\pi_*(\O_{\pa}(p))=0$ et  pour $p=0$,
$R^0\pi_*(\O_{\pa})=\O_{M}$ et $R^q\pi_*(\O_{\pa})=0$ pour $q>0$.
Parmi tous les $R^q\pi_*(\O_{\pa}(p))$ qui apparaissent, le seul non nul est $R^0\pi_*(\O_{\pa})=\O_{M}$. Comme
toutes les flèches sont nulles on a $E_{\infty}=E_3=E_2=E_1$ donc
$R^0\pi_*(\O_{\pf})=\O_{M}$ et $R^q\pi_*(\O_{\pf})=0$
pour $q>0$. D'où i).

Pour ii) on tensorise par $\O(1)$ la suite de Koszul (\ref{equ:(*)}). On
obtient la suite exacte sur $\pa$~:
$$\begin{array}{cccccccccc}
0&\to&\Lambda^{a-r}\pi^*\B(r+1-a)&\to&\Lambda^{a-r-1}\pi^*\B(r+2-a)&\to&\cdots&\to&\Lambda^{2}\pi^*\B(-1)&\to\\
&\to&\pi^*\B&\to&\O_{\pa}(1)&\to&\O_{\pf}(1)&\to&0.&
  \end{array}$$
La restriction de ces \fx aux fibres $\proj^{a-1}$ est respectivement
$\O(r+1-a)^{\cdot}, \ldots, \O(-1)^{\cdot}, \O^{\cdot},\O(1).$ Il
résulte (\cite{Hart}, thm. 5.1) que la cohomologie sur les fibres est
nulle sauf pour $\H^0$ des deux derniers. Par le \th de
semi-continuité il résulte que les images directes supérieures sont
nulles sauf pour $\pi_*(\pi^*\B)=\B=E_1^{-1,0}$ et pour
$\pi_*(\O_{\pa}(1))=\A=E_1^{0,0}$. L'application $d_1^{-1,0}$ entre
ces deux \fx est induite par $\pi_*$ de l'application
$\pi^*\B\to\O_{\pa}(1)$ qui était l'application composée (\ref{equ:AC}):$\pi^*\B\to
\pi^*\A\to\O_{\pa}(1)$. Par $\pi_*$ on a l'égalité
$\pi_*(\O_{\pa}(1))=\A$ (\cite{Hart}, prop. 7.11, p. 162) et le \mor
induit $\pi_*(\pi^*\A)=\A\to \pi_*(\O_{\pa}(1))=\A$ est l'identité. Le
\mor  $\pi_*(\pi^*\B)\to \pi_*(\pi^*\A)$ est le \mor de départ
$\B\to\A$. Au total, le \mor $d_1:\pi_*(\pi^*\B)\to
\pi_*(\O_{\pa}(1))$ est le \mor $\B\to\A$ de départ. Mais ce \mor
satisfait la suite exacte
$$0\to\B\to\A\to\F\to 0.$$
Il résulte que $E_2^{0,0}=\F$ et le reste des termes $E_2$ sont
nuls. Donc 
$$R^q\pi_*\O_{\pf}(1)=\left\{
  \begin{array}{ccc}
\F&\mbox{ si }&q=0\\
0&\mbox{ si }&q>0
  \end{array}
\right..\Box
$$

\begin{cor}
\label{cor:siy}
  Soit $Y$  une sous-variété de Cohen-Macaulay de codimension $2$ de
  $M$ et $\F=I_Y$ le \fs d'idéaux associé. Supposons que l'inclusion de l'éclaté $\blym$ de $M$ selon
  $Y$ dans le schéma en espaces projectifs de Grothendieck $\piy$
  au-dessus de $M$ soit un \iso.
  Alors
 $$R^q\pi_*\O_{E}=\left\{
  \begin{array}{ccc}
\O_Y&\mbox{ si }&q=0\\
0&\mbox{ si }&q>0
  \end{array}
\right.$$
où $E$ est le diviseur exceptionnel.
\end{cor}

\par{\bf Preuve du corollaire~:}
Par la propriété
universelle de $\piy$, comme dans le rappel \ref{pro:soi}, il existe un \mor $\iota:\blym\to\piy$ tel que
$\iota^*(\O(1))=\O(-E)$, et $\iota$ est une immersion fermée, qui est
un \iso par hypothèse. Dans l'\iso $\blym\simeq\piy$,
l'idéal $I_E$ du diviseur exceptionnel $E$ dans $\blym$ est la
préimage du \fs $\O(1)$ sur $\piy$. La suite 
$$0\to I_E\to\O_{\piy}\to\O_E\to 0$$
sur $\piy$ induit une suite longue pour les images directes supérieures
qui implique grâce à la proposition \ref{dcc}~:
$$0\to\pi_*I_E\to\pi_*\O_{\piy}\to\pi_*\O_E\to 0$$
et $R^q\pi_*\O_E=0$ pour $q>0$. Il résulte
$$0\to I_Y\to\O_{M}\to\pi_*\O_E\to 0$$
donc $\pi_*\O_E=\O_{Y}$.$\Box$


\section{La géométrie de la variété d'incidence }
\label{sec:lag}

On fera ici l'étude du \mor $\phi=(p_m,q):\hilxmm\to\hilxx$. On note $\rho:\blup\to\hilxx$ l'éclatement de
$\hilxx$ selon le sous-schéma universel $\xim$ et $\pi:\pix\to\hilxx$
le \mor de projection associé à $\pix$, le schéma en espaces projectifs de Grothendieck
associé au \fs $\ix$. Le but de la section est de
démontrer le
\begin{theor}
\label{the:iso}
On a un diagramme commutatif au-dessus de $\hilxx$~:
$$\diagram
\blup\dto_{}\rto&\pix\dlto^{}_{\simeq}\\
\hilxmm.
\enddiagram$$

\end{theor}

Le \th sera la conclusion des cinq propositions suivantes~:
\begin{defiprop}(\cite{Che},\cite{Tik})
Soit $F_1$ le foncteur
$$F_1:\{\mbox{Schémas}\}\to\{\mbox{Ensembles}\}$$
qui associe à un schéma $S$ l'ensemble des couples $(Z,\xi)$ où
$Z\subset S\times X$ (respectivement $\xi\subset S\times X$) est un
sous-schéma de $S\times X$, $S$-plat, de longueur relative $m+1$
(respectivement $m$) au-dessus de $S$, tels que $\xi\subset Z$. Ce
foncteur est représentable par une variété $\hilxmm$, munie des familles $\xiu\subset \zu$.

\end{defiprop}
 
Ainsi, il existe un schéma $\hilxmm$ et des familles $\zu$ et $\xiu$
au-dessus de $\hilxmm$ telles que 
pour chaque $S$, $Z$, $\xi$, il existe un unique \mor $f:S\to\hilxmm$ vérifiant $(f\times
id)^{-1}(\zu)=Z, (f\times
id)^{-1}(\xiu)=\xi$.

 Par construction, il existe des morphismes  de
projection
$$\diagram
\hilxmm\rto^{p_{m+1}}\dto_{p_{m}}&\hilxp\\
\hilx
\enddiagram$$
qui font de $\hilxmm$ un fermé de $\hilxp\times\hilx$ et
$\zu=(p_{m+1}\times id)^{-1}(\Xi_{m+1})$ et $\xiu=(p_{m}\times id)^{-1}(\Xi_{m})$.

\begin{prop}(cf. \cite{E-S}, lemme 3.1, \cite{Lehn}, pag.8)
  \label{pro:xproj}
Les variétés $\hilxmm$ et $\pix$, le schéma en espaces projectifs de Grothendieck
associé au \fs $\ix$, sont isomorphes au-dessus de $\hilxx$. 
\end{prop}

\par{\bf Preuve de la proposition~:}
On rappelle la propriété universelle de $\pix$~: Pour tout schéma $S$, tout \mor $f:S\to \hilxx$, tout \fib inversible $L$
sur $S$ et tout \mor \sur $f^*\ix\to L\to 0$, il existe un \mor unique
$g:S\to\pix$ qui rend commutatif le diagramme
$$\diagram
&\pix\dto^{\pi}\\
S\urto^{g}\rto_{f\ \ \ \ }&\hilxx
\enddiagram$$
et satisfait $L=g^*\O_{\pix}(1)$. On rappelle aussi que pour une variété $X$, la variété $\hilx$, munie du
sous-schéma universel $\xim$, est l'objet
universel pour la donnée d'un schéma $S$, muni d'une famille $S$-plate
$\xis\subset S\times X$ de sous-schémas finis de
$X$ de longueur $m$  au-dessus de $S$. Sous forme de diagramme cela s'écrit~:
$$
\diagram
 \xis\drto_{(m)}\rto^{\subset}&S\times
 X\dto\\
&S.
\enddiagram $$
Plus précisément, Grothendieck démontre (\cite{Grot2}) qu'il existe une variété $\hilx$ et un
fermé $\xim$ dans $\hilxx$, fini et plat de degré $m$
au-dessus de $\hilx$, tels que pour chaque $(S,\xis)$ comme ci-dessus
il existe un unique \mor $f:S\to\hilx$ vérifiant $(f\times
id)^{-1}(\xim)=\xis$.
Alors la variété $\hilxx$ munie de  $\xipu=(p_1\times
  id)^{-1}(\xim)$ et  $\etapu=(p_2\times
  id)^{-1}(\Xi_1)$ est un objet universel pour la donnée d'un schéma $S$, avec un couple $ (\xis,\etas)$ où
$\xis\subset S\times X$ (respectivement $\etas\subset S\times X$) est
une famille  $S$-plate de sous-schémas finis de
$X$ de longueur $m$ (respectivement $1$) au-dessus de $S$, soit
$$\begin{array}{ccc}
\diagram
 \xis\drto_{(m)}\rto^{\subset}&S\times
 X\dto\\
&S
\enddiagram & \mbox{et}&\diagram
 \etas\drto_{(1)}\rto^{\subset}&S\times
 X\dto\\
&S.
\enddiagram\end{array} $$

Commençons la preuve de la proposition \ref{pro:xproj} par trois lemmes~:
\begin{lemme}
\label{lem:eta}
  Soit $M$  une variété projective et ${\eta}$ un sous-schéma  fermé
  de $S\times M$ qui soit $S$-plat et de longueur relative $1$
  au-dessus de $S$. Le sous-schéma $\eta$ est le graphe d'un \mor  $\etab:S\to M$.
\end{lemme}
\begin{rem}
  {\rm La réciproque est évidente: si $\etab:S\to M$ est un morphisme, le
  graphe $\eta$ de $\etab$  est $S$-plat et de longueur relative $1$ au-dessus de $S$.}
\end{rem}
\par{\bf Preuve du lemme \ref{lem:eta}~:}
L'application $\eta\to S\times M\stackrel{pr_1}{\to}S$ est un \mor
fini et plat de degré $1$, donc un isomorphisme. Son inverse est
de la forme $i:S\to\eta$, $i(s)=(s,\etab(s))$.$\Box$

\begin{lemme}
  \label{lem:sif}
Si $F$ est un \fs cohérent sur $S\times M$, $S$-plat et de longueur
relative $1$ au-dessus de $S$, alors $F$ est de la forme
$\O_{\eta}\tens L$, où $\eta$ est le graphe d'un \mor $\etab:S\to M$
et $L$ un \fs inversible sur $\eta$.
\end{lemme}

\begin{rem}
  \label{rem:eta}
{\rm Dans le lemme  précédent on avait vu que $\eta$ était isomorphe à $S$
par $pr_1$. Par cet \iso on établit une correspondance entre les \fx
sur $S$ et sur $\eta$. On fera des abus de notation en utilisant cette
correspondance. Par exemple si $E$ est un diviseur de Cartier sur $S$,
on note $\O_{\eta}(-E)$ le \fs inversible sur $\eta$ qui correspond à $\O_{S}(-E)$. }
\end{rem}

\par{\bf Preuve du lemme \ref{lem:sif}~:}
Par un résultat de Grothendieck
\cite{EGA}, \S 7, cité par Mumford dans \cite{Mum}, p.19, \S 5 (a)~,
 si $f:X\to Y$ est un \mor propre de schémas noethériens, $\F$ est un
\fs \co sur $X$, plat sur $Y$, pour $y\in Y$, $X_y$ (respectivement
$\F_y$) est la fibre de $f$ au-dessus de $y$ (respectivement le \fs
induit par $\F$ sur la fibre), et si pour tout $y\in Y$,
$\H^1(X_y,\F_y)=0$, alors $f_*(\F)$ est un \fs localement libre sur
$Y$. Dans notre cas,  $A=pr_{1*}F$ est un \fs inversible sur
$S$, parce que $\H^1(M,F_s)=0$ en tout point $s\in S$. Le \mor naturel $pr_1^*A\to F$ est surjectif puisque $F$ est de
longueur relative $1$ au-dessus de $S$, donc le \mor $\O_{S\times
  M}\to (pr_1^*A)^{-1}\tens F$ est surjectif. Le \fs
$(pr_1^*A)^{-1}\tens F$ est encore plat et  de
longueur relative $1$ au-dessus de $S$. On se retrouve dans la
situation du lemme \ref{lem:eta} donc $(pr_1^*A)^{-1}\tens
F=\O_{\eta}$ pour $\eta$ le graphe d'un \mor $\etab:S\to M$. Alors
$F=pr_1^*A\tens\O_{\eta}= L\tens\O_{\eta}$ pour $L=pr_1^*A$, \fs
inversible sur $\eta$. Dans l'abus de notation on peut écrire
$F=\O_{\eta}\tens A$.$\Box$

\begin{lemme}
  \label{lem:lef}
Pour un schéma $S$, la donnée d'un élément de
$F_1(S)$ est équivalente à la donnée de $\xi, \eta, L, s$, où $\xi$ et
$\eta$ sont des sous-schémas $S$-plats de $S\times X$, de longueur
relative $m$, respectivement $1$, au-dessus de $S$, $L$ est un
fibré inversible sur $S$, et $s$ est un \mor \sur $s:I_{\xi}\to L\tens\O_{\eta}$ (dans l'abus de
langage signalé dans la remarque \ref{rem:eta}).
\end{lemme}
\par{\bf Preuve du lemme \ref{lem:lef}~:}
La donnée d'un élément de
$F_1(S)$ nous fournit la suite exacte~:

\begin{equation}
  \label{equ:doi}
0\to I_{\xi}/I_Z\to \O_{S\times X}/I_Z=\O_Z\to \O_{S\times X}/I_{\xi}=\O_{\xi}\to 0.
\end{equation}

La correspondance des données $(\xi, Z)$, avec $\xi\subset Z$ et $(\xi,
\eta, L, s: I_{\xi}\surto L\tens\O_{\eta})$ est obtenue par la suite
exacte (\ref{equ:doi}) en tenant compte du lemme \ref{lem:sif}.$\Box$

Soit $(S,\xi,Z)$ un élément de $F_1(S)$. On note $u:S\to\hilx$ le \mor
donné par la famille $\xi$, $v:S\to X$ le \mor donné par la famille
$\eta$ construite par le lemme \ref{lem:lef}, et
$f=(u,v):S\to\hilxx$. On considère le diagramme commutatif~:
$$\diagram
S\times X\rto^{ u\times id }&\hilxx\\
S.\uto^{\etab=(id,v)}\urto_{f=(u,v)}&
\enddiagram$$
L'image de $\etab=(id,v)$ est $\eta$. Un dernier lemme nous aide à conclure~:
\begin{lemme}
  \label{lem:doi}
On a un \iso de \fx sur $S\times X$, $ (u\times id)^*\ix=I_{\xi}$.
\end{lemme}

\par{\bf Preuve du lemme~:}
À partir de la suite exacte sur $\hilxx$~:
$$0\to\ix\to\O_{\hilxx}\to\O_{\xim}\to 0,$$
par image réciproque par $ u\times id $ on obtient~:
$$\Tor_1^{\O_{\hilxx}}(\O_{S\times X},\O_{\xim})\to
 (u\times id)^*\ix\to\O_{S\times X}\to\O_{\eta}= (u\times id)^*\O_{\xim}\to 0.$$
Puisque $\cdot\tens_{\O_{\hilxx}}\O_{S\times
  X}=\cdot\tens_{\O_{\hilx}}\O_{S}$, il en résulte que
$$\Tor_1^{\O_{\hilxx}}(\O_{S\times
  X},\O_{\xim})=\Tor_1^{\O_{\hilx}}(\O_{S},\O_{\xim})=0,$$ puisque
$\O_{\xim}$ est plat sur $\hilx$. Il en résulte que
$ (u\times id)^*\ix=I_{\xi}$, puisque ce dernier est le noyau de la
surjection $\O_{S\times X}\to\O_{\eta}\to 0$.$\Box$

\par{\bf Retour à la preuve de la proposition \ref{pro:xproj}~:}
D'après le lemme \ref{lem:lef}, se donner un élément de $F_1(S)$
équivaut à se donner une surjection $I_{\xi}\to\O_{\eta}\tens L\to 0$
sur $S\times X$. Mais $\O_{\eta}\tens L$ est un \fs sur le fermé
$\eta=\etab(S)$ de $S\times X$. La surjection
$I_{\xi}\to\O_{\eta}\tens L\to 0$ est équivalente à une surjection
$I_{\xi}\tens_{\O_{S\times X}}\O_{\eta}\to\O_{\eta}\tens L\to 0$. Mais
$\etab$ est un \iso entre $S$ et $\eta$. Par suite la dernière
surjection est équivalente à une surjection $\etab^*I_{\xi}\to
L=\etab^*( \O_{\eta}\tens L)\to 0$. Par le lemme \ref{lem:doi} on a $
\etab^*I_{\xi}=\etab^*(u\times id)^*\ix=f^*\ix.$ Il en résulte que se donner
un élément de $F_1(S)$ revient à se donner un \mor $f:S\to\hilxx$, un
\fs inversible $L$ sur $S$, et une surjection $f^*\ix\to L.$ Par
conséquent $\pix$ représente $F_1$, et donc $\pix\simeq\hilxmm$.$\Box$

Pour démontrer le \th \ref{the:iso} on donnera un \mor entre  les variétés
$\blup$ et $\hilxmm$ (prop. \ref{pro:ile}), on utilisera le morphisme
entre $\blup$ et $\pix$
(rappel \ref{pro:soi}) et on montrera la commutativité du diagramme.
 De plus,
d'après la preuve de la proposition \ref{pro:ile} on obtient sur
$\blup\times X$, pour le diviseur exceptionnel $E$, et pour la famille
$\zupp$ construite au-dessus de $\blup$ à partir des familles $\xippu$
et $\etappu$, la suite
exacte
\begin{equation}
  \label{equ:sir}
  0\to\O_{\zupp}\to\O_{\xippu}\oplus\O_{\etappu}\to \O_{\etappu}|_E\to 0.
\end{equation}

La variété $\blup$ avec le \mor $\rho$
est objet universel pour la donnée d'une variété $S$ et  d'un \mor
$f:S\to \hilxx$ tel que $f^{-1}\ix\cdot\O_S$ est un \fs inversible. 
\begin{prop}
\label{pro:f2}
   La
  variété $\blup$ munie  des familles $\xippu=(\rho\times id)^{-1}((p_1\times
  id)^{-1}(\xim))$ et  $\etappu=(\rho\times id)^{-1}((p_2\times
  id)^{-1}(\Xi_1))$ et du diviseur exceptionnel $E=\rho^{-1}(\xim)$
  est objet universel  pour la donnée d'une variété $S$, des
  familles $S$-plates $\xi\subset S\times X$, $\etas\subset S\times
  X$ de sous-schémas finis de
$X$ de longueur $m$, respectivement $1$, au-dessus $S$ et d'un
  diviseur de Cartier $E_S$, tels que la suite~:
\begin{equation}
  \label{equ:***}
  0\to\O_{\etas}(-E_S)\to\O_{\etas}\to\O_{\eta\cap\xi}\to 0
\end{equation}
est exacte.
\end{prop}

Plus précisément, pour une telle variété $S$, il existe un unique \mor
$f:S\to\blup$, tel que $(f\times
id)^{-1}(\xippu)=\xi, (f\times
id)^{-1}(\etappu)=\etas$ et $f^{-1}(E)=E_S$.

La notation $\O_{\etas}(-E_S)$ est expliquée dans la remarque \ref{rem:eta}.

\par{\bf Preuve de la proposition \ref{pro:f2}~:}

Le \mor \sur $\O_{\hilx\times X\times X}\to\O_{\xipu}$ induit par
l'inclusion de $\xipu$ dans $\hilxx\times X$, tensorisé par
$\O_{\etapu}$, donne le \mor surjectif~:
$$\O_{\etapu}\to\O_{\xipu}\tens_{\O_{\hilxx\times
    X}}\O_{\etapu}=\O_{\xipu\cap\etapu}$$
de noyau $\iu$.
Le \fs $\O_{\xipu}\tens_{\O_{\hilxx\times
    X}}\O_{\etapu}$ est le \fs structural de l'intersection
    $\xipu\cap\etapu$. Le \fs $\iu$ est un \fs d'idéaux dans
    $\O_{\etapu}$. Par la correspondance de la remarque \ref{rem:eta},
    il correspond à un \fs d'idéaux sur $\hilxx$.

    \begin{lemme}
      \label{lem:iu}
Le \fs $\O_{\xipu}\tens_{\O_{\hilxx\times
    X}}\O_{\etapu}$ correspond par l'isomorphisme $\etapu\to\hilxx$ au \fs $\O_{\xim}$. De manière équivalente, l'idéal sur
    $\hilxx$ qui correspond au noyau $\iu$ de la surjection $\O_{\etapu}\to\O_{\xipu\cap\etapu}$, est $\ix$.
    \end{lemme}

\par{\bf Preuve du lemme \ref{lem:iu}~:}
On a, pour $p_{13}, p_{23}, p_{12}$ les projections naturelles, 
$$\begin{array}{c}
  \O_{\xipu}=p_{13}^*\O_{\xim}, \O_{\etapu}=p_{23}^*\O_{\Xi_1}, \\
 \O_{\xipu}\tens_{\O_{\hilxx\times
    X}}\O_{\etapu}=p_{13}^*\O_{\xim}\tens p_{23}^*\O_{\Xi_1}=p_{12}^*\O_{\xim}\tens\O_{\etapu}.
\end{array}$$
Comme $ p_{12}$ est l'\iso entre $\etapu$ et $\hilxx$, le \fs sur
$\hilxx$ qui correspond à $p_{12}^*\O_{\xim}\tens\O_{\etapu}$
est $\O_{\xim}$.$\Box$

\par{\bf Retour à la preuve de la proposition \ref{pro:f2}~:}
Soient $S, \xis, \etas, E_S$ donnés. Par la propriété d'universalité de
$\hilxx$, il existe un unique \mor $f:S\to\hilxx$ tel que
$(f\times id)^{-1}\O_{\xipu}=\O_{\xis}$ et $(f\times
id)^{-1}\O_{\etapu}=\O_{\etas}$. Le \fs d'idéaux correspondant à
$\Ker(\O_{\etas}\to\O_{\eta\cap\xi})$ est
$f_S^{-1}(\Ker(\O_{\etapu}\to\O_{\etapu}\tens
\O_{\xipu}))\cdot\O_S=f_S^{-1}I_{\xim}\cdot\O_S$. Il est de la forme $\O_S(-E_S)$, avec $E_S$
diviseur de Cartier, donc $f$ se factorise à travers l'éclatement
$\rho:\blup\to\hilxx$.$\Box$

\begin{prop}
\label{pro:ile}
  Il existe un \mor canonique $\blup\to\hilxmm$ au-dessus de $\hilxx$.
\end{prop}
\par{\bf Preuve de la proposition~:}

Il suffit de donner une manière de construire  à partir de  $S$, $\xi\subset S\times X$ et $\etas\subset S\times
  X$ des
  familles $S$-plates de sous-schémas finis de
$X$ de longueur $m$, respectivement $1$, au-dessus $S$, et $E_S$ un
  diviseur de Cartier, tels que la suite (\ref{equ:***})
est exacte,  des sous-schémas $Z\subset S\times X$ et $\xi\subset
S\times X$  de
$S\times X$,
$S$-plats, de longueur relative $m+1$, respectivement $m$,
au-dessus de $S$, tels que $\xi\subset Z$. 

Effectivement, par la propriété universelle de $\blup$
(prop. \ref{pro:f2}), on obtient un \mor
$a:\blup\to\hilxmm$ tel que les images réciproques par $a\times id$ des
familles $\zu$, $\xiu$ au-dessus de $\hilxmm$ soient les familles  $\zupp$, $\xippu$ construites sur $\blup$ à partir
de $E$, $\xippu$, $\etappu$.
Décrivons cette méthode. On se donne $S$, $E_S$, $\xi$, $\eta$,
  comme ci-dessus. À partir de la surjection $\O_{S\times
  X}\to\O_{\etas}\to 0$ on obtient la surjection
  $\O_{\xis}\to\O_{\xi\cap\eta}\to 0.$ Soit $\alpha$ le noyau
  du \mor
  \begin{equation}
    \label{equ:*}
    \begin{array}{ccccc}
\O_{\xis}\oplus\O_{\etas}&\to&\O_{\xi\cap\eta}&\to& 0\\
(u,v)&\mapsto&(u-v)|_{\xis\cap\etas}&&.
    \end{array}
  \end{equation}
On dispose également d'un \mor 
$$    \begin{array}{ccc}
\O_{S\times X}&\to&\O_{\xis}\oplus\O_{\etas}\\
t&\mapsto&(t|_{\xis}, t|_{\etas})
\end{array}$$
qui donne $0$ par composition avec le \mor (\ref{equ:*}). D'où un \mor
$\O_{S\times X}\to\alpha$ qui est surjectif. Par conséquent 
\begin{equation}
\label{equ:**}
\alpha=\O_{\zs},
\end{equation}
 où $\zs$ est un sous-schéma fermé de $S\times X$.

On applique le lemme du serpent pour le diagramme
\begin{equation}
\label{dia:ste}
{\diagram
&0\dto&0\dto&0\dto&\\
&\O_{\etas}(-E_S)\dto&\O_{\zs}\dto&\O_{\xis}\dto&\\
0\rto&\O_{\etas}\dto\rto&\O_{\xis}\oplus\O_{\etas}\dto\rto&\O_{\xis}\dto\rto&0\\
0\rto&\O_{\eta\cap\xi}\dto\rto&\O_{\eta\cap\xi}\dto\rto&0\rto&0\\
&0&0&&
\enddiagram}\end{equation}
où la colonne de gauche est donnée par hypothèse. On obtient la suite exacte
$$0\to\O_{\etas}(-E_S)\to\O_{\zs}\to\O_{\xis}\to 0,$$
où $\zs$ est une famille $S$-plate dans $S\times X$, de longueur
relative $m+1$. On a obtenu $Z$ satisfaisant la propriété énoncée. 

La suite exacte (\ref{equ:sir}) résulte pour $S=\blup$ de la suite
(\ref{equ:*}), de l'égalité (\ref{equ:**}) et de la suite exacte
(\ref{equ:***}) donnée dans la proposition \ref{pro:f2}, qui dit
exactement que $\O_{\eta\cap\xi}\simeq\O_{\etas}|_E$.$\Box$

\par{\bf Preuve du \th \ref{the:iso}~:}
Soit $S=\blup$. Dans la démonstration de la proposition \ref{pro:ile}
on avait construit, en partant des données $E, \xi, \eta$, un fermé
$Z$ de $S\times X$, $S$-plat de longueur relative $m+1$ au-dessus de
$S$, tel que la suite 
\begin{equation}
\label{equ:sui}
0\to\O_{\eta}(-E)\to\O_Z\to\O_{\xi}\to 0
\end{equation}
soit exacte.

\begin{lemme}
  \label{lem:lat}
La tensorisation de la suite (\ref{equ:sui}) par $\O_{\eta}$ sur
$\O_{S\times X}$ est~:
$$0\to\O_{\eta}(-E)\to\O_{\eta}\to\O_{\eta\cap\xi}\to 0$$
avec les morphismes canoniques (l'inclusion, respectivement la projection).
\end{lemme}

\par{\bf Preuve du lemme~:}
Dans la démonstration de la proposition \ref{pro:ile} on avait vu que
le \mor $\O_Z\to\O_{\xi}\oplus\O_{\eta}$ provient par restriction à
$Z$ du \mor
$$    \begin{array}{ccc}
\O_{S\times X}&\to&\O_{\xis}\oplus\O_{\etas}\\
t&\mapsto&(t|_{\xis}, t|_{\etas})
\end{array}$$
Il résulte que les morphismes $\O_Z\to\O_{\xi}$ et
$\O_Z\to\O_{\eta}$ qui composent le \mor
$\O_Z\to\O_{\xi}\oplus\O_{\eta}$ du diagramme (\ref{dia:ste}),
proviennent des morphismes canoniques $\O_{S\times
  X}\to\O_Z\to\O_{\xi}$ et $\O_{S\times X}\to\O_Z\to\O_{\eta}$ (ce sont
des restrictions du faisceau structural de $Z$ aux fermés $\xi$, $\eta$). Donc $\O_Z\tens\O_{\eta}=\O_{\eta}$. Ceci
implique, par la commutativité du carré en haut à droite dans le
diagramme (\ref{dia:ste}), que le \mor
$\O_Z\to\O_{\xi}$ de la première ligne est canonique, donc par
tensorisation avec $\O_{\eta}$ c'est le \mor canonique
$\O_{\eta}\to\O_{\eta\cap\xi}$. Si on tensorise par $\O_{\eta}$
le carré en haut à gauche du même diagramme on obtient
$$\diagram
\O_{\eta}(-E)\dto&\O_{\eta}\dto\\
\O_{\eta}\rto&(\O_{\eta\cap\xi})\oplus\O_{\eta}.
\enddiagram$$
Tous les morphismes sont canoniques, donc le \mor
$\O_{\eta}(-E)\to\O_{\eta}$ de la première ligne est l'inclusion
canonique, par la commutativité du diagramme. $\Box$

\par{\bf Retour à la preuve du \th \ref{the:iso}~:}
Le \mor composé $S=\blup\to\hilxmm\to\pix$ est construit explicitement
en quatre pas~:

a) On part de la donnée de $E,\xi, \eta$  et on obtient
$0\to\O_{\eta}(-E)\to\O_Z\to\O_{\xi}\to 0$  dans le diagramme
(\ref{dia:ste}).

b) D'ici on applique la suite (\ref{equ:doi}) et on obtient une
surjection $I_{\xi}\to \O_{\eta}(-E)$.

c) On tensorise avec $\O_{\eta}$ et on obtient une surjection
$I_{\xi}\tens\O_{\eta}\to \O_{\eta}(-E)$.

d) On applique $\etab^*$ à cette surjection et on obtient une
surjection $f^*\ix\to \O_S(-E)$ (les notations et les explications sont
données dans le lemme \ref{lem:doi}: $f$ est le \mor d'éclatement
$f:\blup\to\hilxx$, et $f^*\ix=\etab^*(u\times id)^*\ix= \etab^*I_{\xi}=\etab^*(I_{\xi}\tens\O_{\eta})$).

On remarque que au pas c) le \mor $I_{\xi}\tens\O_{\eta}\to
\O_{\eta}(-E)$ provient par tensorisation par $\O_{\eta}$ de la suite
$0\to I_Z\to I_{\xi}\to \O_{\eta}(-E)\to 0$ déduite de la suite 
(\ref{equ:doi}). D'après le lemme \ref{lem:lat} il résulte que par
tensorisation par $\O_{\eta}$ on obtient un diagramme~:
$$\diagram
&I_{\xi}\tens\O_{\eta}\rto&\O_{\eta}\rto\dto^{\simeq}&\O_{\eta\cap\xi}\dto^{\simeq}\\
0\rto&\O_{\eta}(-E)\rto&\O_{\eta}\rto&\O_{\eta\cap\xi}.
\enddiagram$$
L'image réciproque par $\etab^*$ de ce diagramme est exactement le
diagramme (\ref{equ:unu}) du rappel \ref{pro:soi}. Par suite les morphismes $b:\blup\to\pix$
et  $\blup\stackrel{a}{\to}\hilxmm\to\pix$ sont représentés par la
même surjection $f^*\ix\to\O(-E)$ sur $\blup$, donc ils
coïncident.$\Box$


\section{Images directes du \fs structural  de la variété d'incidence }
\label{sec:lageo}

On rappelle que $\phi$ est le \mor $\phi=(p_m,q):\hilxmm\to\hilxx$. Le
\mor $q$ est rigoureusement défini à partir de la famille $\etau$
obtenue sur $\hilxmm$ des familles $\zu$, $\xiu$ comme dans le lemme
\ref{lem:lef}. Le but de la section est de démontrer la proposition suivante~:

\begin{prop}
  \label{the:3}
On a~:

i) $R^q\phi_*(\O_{\hilxmm})=\left\{
  \begin{array}{ccc}
0&\mbox{ si }&q>0\\
\O_{\hilx\times X}&\mbox{ si }&q=0
  \end{array}\right.$

ii) $R^q\phi_*\O_{E}=\left\{
  \begin{array}{ccc}
0&\mbox{ si }&q>0\\
\O_{\Xi_m}&\mbox{ si }&q=0
  \end{array}\right.$
\end{prop}

Le diviseur exceptionnel $E$ sur $\hilxmm$ s'obtient après
l'identification entre $\blup$ et $\hilxmm$ au-dessus de $\hilxx$. Il
s'agit d'un cas particulier de la proposition \ref{dcc} et du corollaire \ref{cor:siy}. Afin
d'appliquer ce corollaire, on aura besoin d'identifier $\blup$ et 
 $\hilxmm$ au-dessus de $\hilxx$. Ceci est une conséquence de
 l'identification entre  $\blup$ et $\proj(\ix)$ au-dessus de $\hilxx$
 ci-dessous (lemme \ref{the:4}) et du résultat de la section
 précédente (\th \ref{the:iso}, cf. \cite{E-S}, prop. 2.2).  On rappelle que  $\ix$ est l'idéal du sous-schéma $\xim$ dans
$\hilxx$ et $\pix$ le schéma en espaces projectifs de Grothendieck
associé au \fs $\ix$.

Il faut d'abord s'assurer que~:
\begin{lemme}
\label{lem:ab}
Le \fs $\ix$ admet une présentation
$$0\to\B\to\A\to\ix\to 0$$
où $\A$ et $\B$ sont localement libres sur $\hilxx$.
\end{lemme}

\par{\bf Preuve du lemme~:}
Soit $0\to\B\to\A\to\ix\to 0$ une présentation, où $\A$ est localement
libre et $\B$ est le noyau du \mor $\A\to\ix\to 0$. On prouve que $\B$
est localement libre. On note $a$ le rang de $\A$. Soit
$s\in\hilx$. On a $\I_s=\ix|_{\{s\}\times X}$. Comme $\ix$ est plat sur $\hilx$, on obtient la suite exacte~:
$$0\to\B|_{\{s\}\times X}\to\A|_{\{s\}\times X}\to\I_s\to 0.$$
Mais $\I_s$ est de dimension homologique $1$, donc $\B|_{\{s\}\times
  X}$ est localement libre de rang $a-1$. Par conséquent le \fs
cohérent $\B$ a toutes ses fibres de dimension $a-1$, donc il est
localement libre de rang $a-1$.$\Box$

\begin{lemme}(cf. \cite{Lehn}, p.8, \cite{E-S}, p.5)
\label{the:4}
  Les variétés $\blup$ et $\proj(\ix)$ sont isomorphes
  au-dessus de $\hilxx$.
\end{lemme}

\par{\bf Preuve du lemme \ref{the:4}~:}
D'après la section préliminaire \ref{sec:surpf}, dans le cas qui nous intéresse: $\F=\ix$, on obtient que $\pix$ est le
fermé donné par l'annulation de la section $\sigma$ du \fib
$\pi^*\B^*(1)$ de rang $a-1$ sur la variété $\pa$, de dimension
$2m+2+a-1$.

Le lemme de Krull affirme que la dimension de toutes les composantes
irréductibles de $\pix$ est supérieure ou égale à
$(2m+2+a-1)-(a-1)=2m+2$. Il suffit de montrer que $\pix$ est une
variété intègre de dimension $2m+2$.  Alors, par la propriété
universelle de $\pix$, d'après le rappel \ref{pro:soi}, il existe un \mor $\iota:\blup\to\pix$ tel que
$\iota^*(\O(1))=\O(-E)$, et $\iota$ est une immersion fermée. Mais
$\blup$ est de dimension $2m+2$ puisque birationnelle à $\hilxx$. D'où
l'égalité $\blup=\pix$. Ces variétés seront identifiées dans la suite,
ainsi que les morphismes $\phi$ et $\pi$.

On commence par montrer que $\pix$ est une variété irréductible de
dimension $2m+2$. Par le lemme \ref{lic} elle est  localement
intersection complète, ce qui implique en particulier par
\cite{Hart}, prop. 8.23(a) (p. 186) que c'est une variété de
Cohen-Macaulay. On remarque ensuite que l'ensemble de ses points
singuliers est de codimension $\ge 1$. Par le critère de Serre
(\cite{Mats}, p.183), on déduit que $\pix$ est une variété réduite, donc
intègre.

\begin{slemme}(cf. \cite{E-S}, prop.3.2(a))
  La variété $\pix$ est irréductible de
dimension $2m+2$ .
\end{slemme}

\par{\bf Preuve du sous-lemme~:}
On utilise le résultat de la proposition 5 de l'article \cite{E-L},
qu'on applique pour $E=\O_X$, $r=1$, $l=m$, $\mbox{Quot}(E,l)=\hilx$,
$\N=\ix$, $Z=\pix$, $Y_m=\hilxx$ et $Y_{m,i}$ l'ensemble localement
fermé $\{y=(s,x)\in\hilxx| e((\ix)_{s,x})=i+1\}$ avec la structure
réduite. Par définition, si $F$ est un \fs cohérent de $\O_X$-modules,
on note $e(F_x)=\hom_X(F,k(x))$ la dimension de l'espace vectoriel
$F(x)=F\tens_{\O_X}k(x)$, qui par le lemme de Nakayama coïncide avec
le nombre minimal de générateurs de la fibre des germes $F_x$. La
définition de $Y_{m,i}$ est équivalente à celle de la page 4 de
\cite{E-L}, en utilisant le lemme 2 de \cite{E-L}. Comme il est dit à
la page 5, la fibre en $(s,x)\in Y_{m,i}$ du \mor $\pi:\pix\to\hilxx$,
noté $\phi$ par les auteurs, est
$\proj((\ix)_s\tens_{\O_X}k(x))=\proj^i$. La proposition 5 dit que la
codimension de $Y_{m,i}$ dans $Y_m=\hilxx$ est au moins $2i$, soit que
$\dim Y_{m,i}\le 2m+2-2i$. Il résulte que l'adhérence des fibres de
$\pi$ au-dessus de $Y_{m,i}$ sont de dimension au plus
$2m+2-2i+i=2m-i+2$. L'important c'est que pour $i>0$ ces ensembles
sont de dimension strictement inférieure à $2m+2$.

L'ensemble $Y_{m,0}$ est le complémentaire de $\xim$ dans $\hilxx$,
puisque pour $x\not\in\supp(s)$, on a
$(\ix)_s\tens_{\O_X}k(x)=\O_X\tens_{\O_X}k(x)=k(x)$, de dimension
$1$. Inversement, si $\dim((\ix)_s\tens_{\O_X}k(x))=1$, comme
$(\ix)_s\tens_{\O_X}k(x)=((\ix)_s\tens_{\O_X}\O_{X,x})\tens_{\O_{X,x}}k(x)=(\ix)_{s,x}\tens_{\O_{X,x}}k(x)$,
par le lemme de Nakayama, $(\ix)_{s,x}$ est engendré par un générateur
$g$. Si ce générateur est inversible dans $\O_{X,x}$, alors
$(\ix)_{s,x}=(1)$ et $\ix$ est nul au voisinage de $x$ donc
$x\not\in\supp(s)$. Si $g$ s'annule en $x$, alors $(\ix)_{s,x}=(g)$, 
donc localement en $x$, $s$ est donnée par le diviseur défini par
$g=0$, contradiction (puisque $s$ est un schéma de dimension $0$).

Le \mor $\pi$ est un \iso au-dessus de $Y_{m,0}$, puisque
 $\ix|_{Y_{m,0}}$ est trivial.

En résumé: On a un \mor $\pi:\pix\to\hilxx$ et $\hilxx$ est réunion
disjointe des ensembles localement fermés $Y_{m,i}$. Par suite $\pix$
est réunion disjointe des ensembles localement fermés
$\phi^{-1}(Y_{m,i})$. Toutes les parties $\phi^{-1}(Y_{m,i})$ pour
$i>0$ sont de dimension $<2m+2$, et $\phi^{-1}(Y_{m,0})$ est de
dimension $2m+2$ ($\pi$ est un \iso au-dessus de $Y_{m,0}$, et
$Y_{m,0}$ est un ouvert de $\hilxx$). On avait vu que $\pix$ était le
schéma des zéros d'une section dans un \fib, et que par le lemme de
Krull, toutes ses composantes irréductibles étaient de dimension $\ge
2m+2$.

À part l'adhérence de $\phi^{-1}(Y_{m,0})$, qui est irréductible
puisque $Y_{m,0}$ l'est, et qui est de dimension $2m+2$, il ne peut
pas y avoir d'autres composantes irréductibles. (Précisément,
$\pix=\cup_i Y_{m,i}=\cup_i \bar{Y}_{m,i}$. Mais quand un
ensemble est recouvert par un nombre fini de fermés, ses composantes irréductibles
sont les éléments maximaux parmi
ces fermés. S'il existait un $i$ tel que
$\bar{Y}_{m,i}\not\subset\bar{Y}_{m,0}$, $\bar{Y}_{m,i}$ engendrerait
une composante irréductible de dimension $<2m+2$, contradiction. Donc
$\bar{Y}_{m,i}\subset\bar{Y}_{m,0}$, pour tout $i$, d'où
$\pix=\bar{Y}_{m,0}$.)$\Box$ 

On remarque que  $Y_{m,0}$ est un ouvert de
$\hilxx$, donc il est lisse, donc l'ensemble des points
singuliers de $\pix$ est de codimension $\ge 1$.$\Box$


\section{Le \mor trace}
\label{sec:lem}

On utilisera ici la propriété de lissité de $\hilxmm$, (voir
\cite{Che}, \cite{Tik}), annoncée depuis l'introduction.
Le but de cette section est de démontrer~:
\begin{prop}
\label{pro:tr}
  L'espace de \coh $\H^q(\hilxp,\lmp\tens\dep)$ est facteur direct
  de \linebreak
  $\H^q(\hilxmm,p_{m+1}^*(\lmp\tens\dep))$. 
\end{prop}

C'est un corollaire de la proposition suivante~:
\begin{prop}
\label{pro:fdir}
  Soient $X$ et $Y$ deux variétés analytiques compactes de même
  dimension $d$, lisses et irréductibles. Soit $f:X\to Y$ un \mor \sur et
  $V$ un \fib vectoriel analytique sur $Y$. Alors pour tout $q\ge 0$,
  $\H^q(Y,V)$ est facteur direct dans $\H^q(X,f^*V)$.
\end{prop}

\begin{cor}
  Si $X$ et $Y$ sont des variétés projectives et lisses de même
  dimension, $f:X\to Y$ un \mor \sur et
  $V$ un \fib vectoriel \alg sur $Y$, alors $\H^q(Y,V)$ est facteur
  direct dans $\H^q(X,f^*V)$ pour tout $q\ge 0$.
\end{cor}

\par{\bf Preuve du corollaire~:}
Puisque $Y$ est projective et $V$ est un \fib algébrique, le principe GAGA
s'applique, et la \coh  $\H^q(Y,V)$ du \fib \alg $V$ sur $Y$ coïncide avec la \coh de $V$ vu comme \fib analytique sur la variété
analytique $Y$. De même pour $X$, $f^*V$. La proposition
\ref{pro:fdir}  s'applique.$\Box$

\par{\bf Preuve de la proposition \ref{pro:fdir}~:}
On procédera en plusieurs étapes~:

1) L'espace $\H^q(Y,V)$ est la \coh en degré $q$ du complexe de
Dolbeault des formes différentielles de type $(0,p)$ sur $Y$ à valeurs
dans $V$: $(A^{0,\cdot}(V),\bar{\partial})$.

2) Il existe un \mor bien défini $f^*:\H^q(Y,V)\to\H^q(X,f^*V)$ qui
associe à une forme fermée sur $Y$ sa préimage sur $X$. On dispose de
la dualité de Serre~: le \mor bilinéaire
$\H^q(Y,V)\tens\H^{d-q}(Y,V^*\tens\oy)\to\H^d(Y,\oy)\simeq\comp$ est
non-dégénéré, c'est-à-dire qu'il identifie chacun des espaces
$\H^q(Y,V)$ et $\H^{d-q}(Y,V^*\tens\oy)$ au dual de l'autre.
(Concrètement ce \mor se réalise ainsi~: on représente l'élément de
$\H^q(Y,V)$ par une forme différentielle $\bar{\partial}$-fermée
$\alpha$ de type $(0,q)$ à valeurs dans $V$ -qui s'écrit localement
$\sum_is_i\tens\omega_i$, $s_i$ section locale holomorphe de $V$ et
$\omega_i$ forme différentielle $C^{\infty}$ de type $(0,q)$- et
l'élément de $\H^{d-q}(Y,V^*\tens\oy)$ par une forme différentielle $\bar{\partial}$-fermée
$\beta$ de type $(0,d-q)$ à valeurs dans $V^*\tens\oy$ -qui s'écrit localement
$\sum_is_i\tens\eta_i\tens\omega_i$, $s_i$ section locale holomorphe
de $V^*$, $\eta_i$   section locale holomorphe
de $\oy$ et 
$\omega_i$ forme différentielle $C^{\infty}$ de type $(0,d-q)$-. Si on
regarde les formes différentielles holomorphes comme des formes
$C^{\infty}$, on obtient que $\alpha\beta$ est une forme \diff $(d,d)$
à valeurs dans $V^*\tens V$, soit un élément dans
$\Gamma(V^*\tens_{\O_Y}V\tens_{\O_Y}\A^{d,d}(Y))$. Il existe un \mor
trace $tr:\Hom_{\O_Y}(V,V)\to\O_Y$. L'image $tr(\alpha\beta)$ est un
élément de $A^{d,d}(Y)=\Gamma(\O_Y\tens_{\O_Y}\A^{d,d}(Y))$. L'\app
qui donne la dualité de Serre est
$(\alpha,\beta)\mapsto\int_Y\alpha\beta$.)

3) On note $\oxy$ le \fib \inve sur $X$ égal à $\ox\tens
f^*(\oy)^{-1}$.
\begin{lemme}
  Le \fib $\oxy$ admet une \sct canonique $\sigma$.
\end{lemme}

\par{\bf Preuve du lemme~:}
L'existence d'une section canonique équivaut à l'existence d'un \mor
injectif  $\O_X\hookrightarrow\ox\tens
f^*(\oy)^{-1}$, ou bien, en tensorisant par $f^*(\oy)$, à l'existence
d'un \mor $f^*(\oy)\to\ox$. Ce \mor s'obtient à partir du \mor entre
les espaces tangents $T_X\to f^*T_Y$ en prenant le dual et la
puissance extérieure
maximale~:$f^*\oy=f^*\Lambda^dT_Y^*\to\ox=\Lambda^dT_X^*$. 

4) La dualité de Serre induit un \mor
$$f_*:\H^q(X,f^*V\tens\oxy)\to\H^q(Y,V).$$
Plus précisément
$$
\H^q(X,f^*V\tens\oxy)\stackrel{Serre}{\simeq}\H^{d-q}(X,f^*V^*\tens\oxy^*\tens\ox)^*\simeq\H^{d-q}(X,f^*V^*\tens
f^*\oy)^*$$
et
$$\H^q(Y,V)\stackrel{Serre}{\simeq}\H^{d-q}(Y,V^*\tens\oy)^*.$$
On dispose d'un \mor image réciproque 
$$f^*:\H^{d-q}(Y,V^*\tens\oy)\to\H^{d-q}(X,f^*(V^*\tens\oy)).$$
Compte-tenu des isomorphismes ci-dessus, on trouve par dualité un \mor
$$f_*:\H^q(X,f^*V\tens\oxy)\to\H^q(Y,V).$$
Concrètement, ce \mor associe à la classe $[\gamma]$ représentée par
la forme \diff $(0,q)$ à valeurs dans $f^*V\tens\oxy$, l'unique classe
de \coh de $\H^q(Y,V)$ qui vérifie~: pour toute forme \diff $\delta$
de cette classe, et pour toute $(0,d-q)$-forme \diff
$\bar{\partial}$-fermée $\epsilon$ à valeurs dans $V^*\tens\oy$, on a 
$$\int_X\gamma\cdot f^*\epsilon=\int_Y\delta\cdot\epsilon.$$
La dualité de Serre assure l'existence de cette classe
$[\delta]=f_*[\gamma]$.

5) En composant les \mors obtenus aux étapes $1$, $3$, $4$, on obtient
un \mor
$$\H^q(Y,V)\stackrel{(1),f^*}{\to}\H^q(X,f^*V)\stackrel{(3)}{\to}\H^q(X,f^*V\tens\oxy)\stackrel{(4),f_*}{\to}\H^q(Y,V).$$
D'après le \th de Sard, l'ensemble des valeurs critiques d'un
morphisme propre  de variétés analytiques lisses est un fermé analytique
sans point intérieur
(ce qui signifie pour nous \iso local, les deux variétés étant de même dimension).
On montre que la composée des trois \apps est la multiplication par
une constante $c$ égale au nombre des points de la fibre au-dessus d'un
point lisse. 

Si $[\alpha]$ est une classe de \coh de $\H^q(Y,V)$, il faut démontrer
que $f_*(f^*[\alpha]\cdot[\sigma])=c\cdot[\alpha]$ où $\sigma$ était
la section canonique construite à l'étape ($3$). Le lemme suivant fait
une première simplification~:
\begin{lemme}[formule de projection]
On a $f_*(f^*[\alpha]\cdot[\sigma])=[\alpha]\cdot f_*([\sigma]).$
\end{lemme}
\begin{rem}
  {\rm On voit ici $[\sigma]$ comme élément de $\H^0(X,\oxy)$ qui est
  envoyé comme à l'étape (4) par $f_*$ dans $\H^0(Y,\O_Y)$. Par la
  compacité de $Y$, ce dernier espace est isomorphe à $\comp$.}
\end{rem}

\par{\bf Preuve de la formule de projection~:}

La classe $f_*([\sigma])$ est représenté par une $0$-forme \diff \hol
$\rho$ sur $Y$ (donc une constante) qui vérifie
$$\int_Y\gamma\cdot\epsilon=\int_X\sigma\cdot f^*\epsilon$$
pour tout $\epsilon\in A^{0,d}(Y,\oy)$. Mais pour
$\epsilon=\alpha\cdot\epsilon'$ avec $\epsilon'\in
A^{0,d-q}(Y,V^*\tens\oy)$ on obtient~:
$$\int_Y\alpha\cdot\gamma\cdot\epsilon'=\int_Y\gamma\cdot\alpha\epsilon'=\int_X
\sigma\cdot f^*(\alpha\epsilon')=\int_Xf^*\alpha\cdot\sigma\cdot
f^*\epsilon'$$
c'est-à-dire exactement $$f_*(f^*[\alpha]\cdot[\sigma])=[\alpha]\cdot
f_*([\sigma]). \Box$$ 

Il reste à montrer que $f_*([\sigma])=c\in\H^0(Y,\O_Y)\simeq\comp$.
On sait que $f_*([\sigma])$ est une constante
$cst\in\H^0(Y,\O_Y)\simeq\comp$. Il faut déterminer $cst$. Pour cela
il suffit de choisir $\epsilon\in A^{0,d}(Y,\oy)$ tel que
$\int_Y\epsilon\ne 0$. Alors $f_*([\sigma])=cst$ implique
$$cst\cdot\int_Y\epsilon=\int_Ycst\cdot\epsilon\int_X
\sigma\cdot f^*\epsilon$$
et on peut déterminer $cst$. Plus précisément, on choisit $y\in Y$ un
point de l'ouvert de lissité. Sur un ouvert dans la topologie
$C^{\infty}$ au-dessus de $y$, $f$ est la réunion de $c$ isomorphismes
$W_i\stackrel{f}{\to}T$, $f^{-1}(T)=\cup^c_{i=1}W_i$.
On choisit un système de coordonnées $z_1,\ldots,z_d$ sur $T$ qu'on
transporte par $f^{-1}$ sur chacun des $W_i$. On choisit, dans cet
abus de notations pour les coordonnées, $dz_1\cdots dz_d$ un
générateur pour $\oy|_T$ et pour $\ox|_{W_i}$. Alors $\sigma=1$ dans
cette trivialisation. Pour $\epsilon$ on choisit une forme 
à support compact inclus dans $T$ et telle que
$\int_T\epsilon=1$. Alors
$$\int_X f^*\epsilon\cdot\sigma=\int_{f^{-1}(T)}
f^*\epsilon\cdot\sigma=\sum^c_{i=1}\int_{W_i}f^*\epsilon=$$
$$=\sum^c_{i=1}\int_T\epsilon=c\cdot\int_T\epsilon.$$
Comme $\int_Y\epsilon=1\ne 0$ alors $cst=c$ et $c\ne 0$ puisque par
hypothèse, $f$ était surjectif. Ceci achève la démonstration de la proposition
\ref{pro:fdir}. $\Box$

Une variante \alg particulière de la proposition qu'on vient de
démontrer est la suivante~:
\begin{prop}
  \label{pro:tra}
Soit $f:X\to Y$ un \mor fini, plat, de degré $c$, entre deux schémas
$X$ et $Y$. Soit $L$ un \fs localement libre sur $Y$. Alors
$\H^q(Y,L)$ est facteur direct de $\H^q(Y,L\tens f_*\O_X)$. 
\end{prop}

Cette proposition est une conséquence du lemme~:
\begin{lemme}
  Soit $f:X\to Y$ comme dans la proposition. Alors $\O_Y$ est facteur
  direct de $f_*\O_X$.
\end{lemme}

Effectivement, soit 
$$\O_Y\stackrel{i}{\hookrightarrow}f_*\O_X\stackrel{p}{\to}
\O_Y$$
 les \mors d'inclusion et projection donnés par le lemme, avec
$p\circ i=id$. Alors la composée des \apps induites 
$$\H^q(Y,L)\to\H^q(Y,L\tens f_*\O_X)\to\H^q(Y,L)$$
est égale à l'identité, d'où le résultat de la proposition.

\par{\bf Preuve du lemme~:} 
Le \fs $f_*\O_X$ est localement libre de rang $c$ sur $Y$. La
multiplication induit un \mor
$$f_*\O_X\to\underline{\mbox{End}}_Y(f_*\O_X,f_*\O_X),$$
qui par composition avec le \mor trace
$$Tr:\underline{\mbox{End}}_Y(f_*\O_X,f_*\O_X)\to \O_Y$$
donne le \mor $p:f_*\O_X\to\O_Y$. Le \mor composé
$$f:\O_Y\stackrel{i}{\hookrightarrow}f_*\O_X\stackrel{p}{\to}
\O_Y$$
est égal à $c\cdot id$. Modulo une reparamétrisation on obtient le résultat. $\Box$


\section{Preuve du \th \ref{the:1}}
\label{sec:preu}

On raisonne par r{\'e}currence sur $m$. On {\'e}tablira tout d'abord les 6
pas suivants~:

1) On a une suite exacte (\ref{equ:sir}) sur $\hilxmm\times X$~:
$$ 0\to\O_{\zu}\to\O_{\xiu}\oplus\O_{\etau}\to \O_{\etau}|_E\to 0$$
o{\`u} $\O_{\etau}|_E$ est le \fs correspondant {\`a} $\O_E$ sur $\hilxmm$ par
l'\iso $\etau\to\hilxmm$.

2)  On a une suite exacte sur $\hilxmm$~:
$$0\to p_{m+1}^*(\lmp\tens\dep)\to \phi^*\left((\lm\tens\de)\boxtimes
A\right)\oplus\phi^*\left(\de\boxtimes (L\tens A)\right)\to$$
$$\to \phi^*\left(\de\boxtimes (L\tens
A)\right)|_E\to 0$$
o{\`u} $\phi=(p_m,q)$.

3) On a une suite longue de \coh associ{\'e}e {\`a} la suite exacte de (2).

4) On peut ramener le calcul de \coh sur $\hilxmm$ {\`a} un calcul de \coh
 sur $\hilxx$~:
$$\H^q(\hilxmm,\phi^*\left((\lm\tens\de)\boxtimes
A\right)\oplus\phi^*\left(\de\boxtimes (L\tens A)\right))=$$
$$=\H^q(\hilxx,(\lm\tens\de)\boxtimes A)\oplus\H^q(\hilxx,\de\boxtimes
(L\tens A)).$$

5) L'espace $\H^q(\hilx,\lmp\tens\dep) $ est facteur direct de
$\H^q(\hilxmm, p_{m+1}^*(\lmp\tens\dep))$.

6) On a
$$\H^q(\hilxmm,\phi^*\left(\de\boxtimes (L\tens
A)\right)|_E)=\H^q(E,\phi^*(\de\boxtimes
(L\tens A)))=$$
$$=\H^q(\xim, p_1^*\de\tens p_2^*(L\tens
A))=\H^q(\hilx,\lam\tens\de).$$

La suite du premier pas est la suite (\ref{equ:sir}) obtenue dans la
 section \ref{sec:lag}. Afin d'{\'e}tablir le deuxi{\`e}me pas, on aura besoin
 de quelques notations et r{\'e}sultats pr{\'e}liminaires.

On s'{\'e}tait donn{\'e} un \fs \inve $A$ sur $X$. On a une suite d'\apps
$$\xmp\stackrel{\alpha}{\to}\smx\times X\stackrel{\beta}{\to}\smxp$$
et un \fs \inve $A\boxtimes A\boxtimes \cdots\boxtimes A$ sur $\xmp$
muni d'une action de $\sigmp$. On a une \app $\gamma:\xm\to\smx$ et un
\fs \inve $A\boxtimes A\boxtimes \cdots\boxtimes A$ sur $\xm$, muni
d'une action de $\sigm $. En appliquant le lemme de descente, par passage au
quotient du \fib $A\boxtimes A\boxtimes \cdots\boxtimes A$ sur $\xm$, on
obtient un \fs $\dm$ sur $\smx$, v{\'e}rifiant $\gamma^*\dm=A\boxtimes
A\boxtimes \cdots\boxtimes A$ ($m$ fois).

\begin{prop}
\label{pro:ilex}
  Il existe un \fs $\dmp$ sur $\smxp$,
  v{\'e}rifiant $\beta^*\dmp=\dm\boxtimes A$.
\end{prop}
\begin{rem}
  {\rm Par cons{\'e}quent, on a 
$$\alpha^*\beta^*\dmp=(\gamma\times id)^*(\dm\boxtimes
 A)=\gamma^*\dm\boxtimes A=A\boxtimes A\boxtimes \cdots\boxtimes A.$$
}
\end{rem}
On introduit alors un \fib $\dem$ sur $\hilxmm$ par image r{\'e}ciproque
de la mani{\`e}re suivante~: $\dem=\mu^*(\dm\boxtimes
A)$, o{\`u} $\mu$ est d{\'e}fini comme dans  le diagramme
commutatif 
$$\diagram
\hilxp\dto_{HC}&\hilxmm\lto_{p_{m+1}}\dto_{\mu}\rto^{\phi=(p_m,q)}&\hilxx\dto^{HC\times
  id}\\
\smxp&\smx\times X\lto_{\beta}\rto^{id}&\smx\times X.
\enddiagram$$
On obtient
\begin{equation}
\label{equ:dem}
\dem=\mu^*\beta^*\dmp=p_{m+1}^*HC^*\dmp=p_{m+1}^*\dep.
\end{equation}

\par{\bf Preuve de la proposition \ref{pro:ilex}~:}
Le \fib $\dmp$ est d{\'e}fini naturellement par passage au
quotient du \fib $A\boxtimes A\boxtimes \cdots\boxtimes A$ sur $\xmp$. Il faut montrer que $\beta^*\dmp=\dm\boxtimes
A$. Puisque leurs images r{\'e}ciproques dans $\Pic(\xmp)$ sont {\'e}gales, il
suffit de d{\'e}montrer que l'\app $a:\Pic(Y/G)\to\Pic(Y)$ est injective, o{\`u}
$Y=\xmp$, $G=\sigm\times \{id\}$, et $Y/G=\smx\times X$ est le quotient de
$Y$ par l'action de $G$. Notons $\Pic^G(Y)$ le groupe de Picard des
$G$-\fibs inversibles sur $Y$, $\Hom(G,\comp^*)$ le groupe des
caract{\`e}res de $G$, $G_y$ le stabilisateur d'un point $y\in Y$, et  $\Hom(G_y,\comp^*)$ le groupe des
caract{\`e}res de $G_y$.
Dans le diagramme
$$\diagram
&&\prod_{y\in Y}\Hom(G_y,\comp^*)&\\
0\rto&\Hom(G,\comp^*)\rto^{e}&\Pic^G(Y)\uto^{d}\rto^{c}&\Pic(Y)\\
&&\Pic(Y/G)\uto^{b}\urto_{a}&\\
&&0\uto&
\enddiagram$$
l'\app $b$ est l'image r{\'e}ciproque, l'\app $c$ est l'oubli, l'\app
$a$ est la  compos{\'e}e $c\circ b$, l'\app $e$ associe {\`a} chaque caract{\`e}re
de $G$ la $G$-action sur $\O_Y$ induite par la $G$-action sur $Y$
associ{\'e}e, et l'\app $d$ associe {\`a} chaque $G$-\fib $L$ et chaque point
$y\in Y$, le caract{\`e}re associ{\'e}
{\`a} l'action induite par $G_y$ sur la fibre au-dessus de $y$, $L_y$.
Il est facile de voir et il est d{\'e}montr{\'e} dans \cite{LeP-D} que la
suite horizontale est exacte.
L'exactitude de la suite verticale est une cons{\'e}quence du lemme de
Kempf (\cite{LeP-D}) qui dit exactement qu'un $G$-\fib L sur $Y$ provient du
quotient $Y/G$ si et seulement si l'action du stabilisateur $G_y$ sur
la fibre $L_y$ est triviale pour chaque $y\in Y$.
L'\app $b$ est injective, puisque si l'image r{\'e}ciproque du \fib 
F sur $Y/G$ est un $G$-\fib trivial $b^*F=\O_Y$ sur $Y$, elle est
munie d'une section $G$-invariante partout non-nulle, qui descend en
une section partout non-nulle de $F$, qui trivialise $F$.
Pour se convaincre que $a$ est injective, il suffit de v{\'e}rifier que
l'intersection $\Hom(G,\comp^*)\cap \Pic(Y/G)$ dans $\Pic^G(Y)$ est
r{\'e}duite {\`a} l'{\'e}l{\'e}ment trivial avec $G$-action triviale. Mais si $f\in
\Hom(G,\comp^*)$ est telle que $d(e(f))=1$ dans $\prod_{y\in
  Y}\Hom(G_y,\comp^*)$, alors $f|_{G_y}=1$ pour tout $y\in Y$. Si $y$
est un point de la forme $(x_l)_l$ avec $x_i=x_j$ alors la
transposition $\tau_{i,j}\in G$ appartient {\`a} $G_y$, donc
$f(\tau_{i,j})=1$. Mais $G$ est engendr{\'e} par toutes ses
transpositions, donc forc{\'e}ment $f$ est le caract{\`e}re trivial. $\Box$

\begin{prop}
  On peut faire un changement de base dans le diagramme
$$\diagram
\zu\dto_{p_1}\rto^{(p_{m+1}\times
  id)}&\ximp\dto^{p_1}\\
\hilxmm\rto^{p_{m+1}}&\hilxp
\enddiagram$$
pour le \fs $\F=\O_{\ximp}\tens p_2^*L$ sur $\ximp$.
\end{prop}

C'est-{\`a}-dire qu'on a
\begin{equation} 
\label{equ:a}
p_{1*}(\O_{\zu}\tens p_2^*L)=p_{m+1}^*(\lmp).
\end{equation}
Ici, on note $p_2$ la deuxi{\`e}me projection $\hilxp\times X\to X$, aussi
bien que la deuxi{\`e}me projection $\hilxmm\times X\to X$.

Cette proposition est une cons{\'e}quence d'un r{\'e}sultat de Grothendieck
\cite{EGA}, \S 7, cit{\'e} par Mumford dans \cite{Mum}, p.19, \S 5 (a)~:

Si $f:X\to Y$ est un \mor propre de sch{\'e}mas noeth{\'e}riens, $\F$ est un
\fs \co sur $X$, plat sur $Y$, pour $y\in Y$, $X_y$ (respectivement
$\F_y$) est la fibre de $f$ au-dessus de $y$ (respectivement le \fs
induit par $\F$ sur la fibre), et que pour tout $y\in Y$,
$\H^1(X_y,\F_y)=0$, alors $f_*(\F)$ est un \fs localement libre sur
$Y$, et ``la formation de $f_*$ commute avec le changement de base'',
i.e. dans toutes les situations de produit \fib~:
$$\diagram
X'\dto_{f'}\rto^{g'}&X\dto^{f}\\
Y'\rto_{g}&Y
\enddiagram$$
le \mor naturel~:
$$g^*(f_*\F)\to f'_*(g^{\prime *}(\F))$$
est un isomorphisme. Ce r{\'e}sultat d{\'e}coule du \th 7.7.5 et la
proposition 7.7.10 appliqu{\'e}e {\`a} $p=1$, et la proposition
7.8.4. appliqu{\'e}e {\`a} $p=0$ (toutes dans \cite{EGA}). On avait d{\'e}j{\`a}
utilis{\'e} la premi{\`e}re partie de ce r{\'e}sultat dans la d{\'e}monstration du
lemme \ref{lem:sif}.

Toutes ces conditions sont satisfaites dans notre cas, donc la
proposition en d{\'e}coule.

\begin{rem}
  {\rm Le m{\^e}me r{\'e}sultat vaut pour $p_{m+1}$ remplac{\'e} par $p_m$,
  respectivement $q$, $\hilxp$ remplac{\'e} par $\hilx$, respectivement
  $X$, et $\F$ remplac{\'e} par $\O_{\xim}\tens p_2^*L$, respectivement
  $\O_{\Delta}\tens p_2^*L$. On obtient ainsi sur $\hilxmm$~:
\begin{eqnarray}
p_{1*}(\O_{\xiu}\tens p_2^*L)\simeq p_m^*(\lm), \label{equ:b}\\
p_{1*}(\O_{\etau}\tens p_2^*L)\simeq q^*(L)\label{equ:c}
\end{eqnarray}
et ce dernier \iso nous donne
\begin{equation}
\label{equ:d}
p_{1*}(\O_{\etau}|_E\tens p_2^*L)\simeq q^*(L)|_E.
\end{equation}}
\end{rem}

En effet, $p_1$ est un \iso entre $\etau$ et $\hilxmm$ et $p_{1*}$
{\'e}tablit un \iso entre les \fx sur $\etau$ et $\hilxmm$. Le dernier
\iso nous dit que, par l'\iso $p_{1*}$, $\O_{\etau}\tens p_2^*L$
correspond {\`a} $q^*L$. Mais $\O_{\etau}|_E$ correspond {\`a} $\O_E$, donc
$\O_{\etau}|_E\tens p_2^*L$ correspond {\`a} $q^*(\lm)|_E.$

Maintenant, le deuxi{\`e}me pas r{\'e}sulte du premier pas par image directe
par $p_1$ de la suite (\ref{equ:sir}) tensoris{\'e}e par $p_2^*L\tens
p_1^*\dem$. Plus exactement il faut montrer~:

a) $p_{1*}(\O_{\zu}\tens p_2^*L\tens
p_1^*\dem)\simeq p_{m+1}^*(\lmp\tens
\dep);$

b) $p_{1*}(\O_{\xiu}\tens p_2^*L\tens
p_1^*\dem)\simeq p_{m}^*(\lm\tens
\de)\tens q^*A;$

c) $p_{1*}(\O_{\etau}\tens p_2^*L\tens
p_1^*\dem)\simeq p_{m}^*(
\de)\tens q^*(L\tens A);$

d) $p_{1*}(\O_{\etau}|_E\tens p_2^*L\tens
p_1^*\dem)\simeq p_{m}^*(
\de)\tens q^*(L\tens A)|_E.$

Ces formules d{\'e}coulent de celles d{\'e}j{\`a} {\'e}tablies (\ref{equ:a}),
(\ref{equ:b}), (\ref{equ:c}), (\ref{equ:d}), en utilisant la formule
de projection et la description de $\dem$, (\ref{equ:dem}).

On remarque ensuite que
\begin{eqnarray*}
 p_{m}^*(\lm\tens
\de)\tens q^*(L\tens A)&=&\phi^*((\lm\tens
\de)\boxtimes A)\\
p_{m}^*(
\de)\tens q^*(L\tens A)&=& \phi^*(\de\boxtimes (L\tens A))
\end{eqnarray*}
o{\`u} $\phi=(p_m,q)$.

On passe au quatri{\`e}me pas.
Dans la section \ref{sec:lageo} on avait vu que 
$$R^q\phi_*\O_{\hilxmm}=\left\{
  \begin{array}{ccc}
0&\mbox{ si }&q>0\\
\O_{\hilx\times X}&\mbox{ si }&q=0.
  \end{array}\right.$$
Alors 
$$R^q\phi_*(\phi^*((\lm\tens
\de)\boxtimes A))=\left\{
  \begin{array}{ccc}
0&\mbox{ si }&q>0\\
(\lm\tens
\de)\boxtimes A&\mbox{ si }&q=0.
  \end{array}\right.$$
Par la suite spectrale de Leray on obtient
$$\H^q(\hilxmm,\phi^*((\lm\tens
\de)\boxtimes A))=\H^q(\hilxx,(\lm\tens
\de)\boxtimes A).$$
De m{\^e}me
$$\H^q(\hilxmm,\phi^*(\de\boxtimes (L\tens
A)))=\H^q(\hilxx,\de\boxtimes (L\tens A)).$$

Le cinqui{\`e}me pas est une cons{\'e}quence de la section \ref{sec:lem}.

Le sixi{\`e}me pas s'obtient comme le quatri{\`e}me en utilisant le r{\'e}sultat
de la section \ref{sec:lageo}~:
$$R^q\phi_*\O_{E}=\left\{
  \begin{array}{ccc}
0&\mbox{ si }&q>0\\
\O_{\Xi_m}&\mbox{ si }&q=0
  \end{array}\right.$$
Alors
$$R^q\phi_*(\phi^*(\de\boxtimes (L\tens
A))|_E)=\left\{
  \begin{array}{ccc}
0&\mbox{ si }&q>0\\
\de\boxtimes (L\tens
A)|_{\xim}&\mbox{ si }&q=0.
  \end{array}\right.$$
Ainsi
$$\H^q(E,\phi^*(\de\boxtimes (L\tens
A)))=\H^q(\xim, \de\boxtimes (L\tens
A)).$$
Le \mor $\xim\stackrel{p_1}{\to}\hilx$ est fini, donc ses images
directes sup{\'e}rieures sont nulles. En appliquant encore une fois la
suite spectrale de Leray on obtient
$$\H^q(\xim, \de\boxtimes (L\tens
A))=\H^q(\hilx,\de\boxtimes \lam).$$

\par{\bf Preuve du th{\'e}or{\`e}me \ref{the:1}~:}

On utilisera une r{\'e}currence sur $m$, pour tous les $m$ et $L$ {\`a} la
fois, pour $A$ fix{\'e}.

Le cas $m=1$ d{\'e}coule du \th d'annulation de Kodaira. Les cas $q=0$ et
$q=1$ sont d{\'e}montr{\'e}s dans les pr{\'e}liminaires. On suppose connus les cas
$1, m$ et on d{\'e}montre le r{\'e}sultat pour $m+1$, $m\ge 1$.

Le r{\'e}sultat du pas (4) s'{\'e}crit apr{\`e}s application de la formule de
K{\"u}nneth~:
$$\H^q(\hilxmm, p_{m}^*(\lm\tens\de)\tens
q^*A\oplus p_n^*\de\tens q^*(L\tens A))=$$
$$=\oplus_{i+j=q}(\H^i(\hilx, \lm\tens\de)\tens \H^j(X,
A))\oplus\oplus_{i+j=q}(\H^i(\hilx, \de)\tens \H^j(X, L\tens A)).$$
Si $q>0$ on a~: 

ou bien $i>0$ et alors 

\hskip 1cm $\H^i(\hilx, \lm\tens\de)=0$
par hypoth{\`e}se de r{\'e}currence pour $m, L$ (on a $\ox^{-1}\otimes
L\otimes A^{\tens k}$  ample

\hskip 1cm  pour $1\le k\le m+1$ donc aussi pour $1\le k\le
m$, et de m{\^e}me pour $\ox^{-1}\otimes A^{\tens k}$) et

\hskip 1cm $\H^i(\hilx,\de)=0$ puisque par hypoth{\`e}se de r{\'e}currence
pour $m, L=\O$ on a 

\hskip 1cm $\H^i(\hilx,\de\tens\om)= \H^i(\hilx,\de\tens
p_{1*}(\O_{\xim}))=0$ et $\H^i(\hilx,\de)$ est facteur direct dans ce

\hskip 1cm dernier espace de \coh par la proposition \ref{pro:tra} de la section
\ref{sec:lem};

ou bien $j>0$ et alors 

\hskip 1cm  $\H^j(X,A)=0$ par Kodaira puisque
  $\ox^{-1}\tens A$ est ample et 

 \hskip 1cm $\H^j(X,L\tens A)=0$ par Kodaira puisque
 $\ox^{-1}\tens L\tens A$ est ample.

Si $q>0$ on a aussi $\H^q(\hilx,\lam\tens\de)=0$ par hypoth{\`e}se de r{\'e}currence
pour $m, L\tens A$ (on v{\'e}rifie que les conditions d'amplitude
v{\'e}rifi{\'e}es par $\ox^{-1}\tens L\tens A^{\tens k}$ pour $1\le k\le m+1$  sont
v{\'e}rifi{\'e}es par $\ox^{-1}\tens (L\tens A)\tens A^{\tens k}$ pour $1\le k\le m$).

Par cons{\'e}quent, dans la suite longue de \coh du pas (3), on obtient
$$\H^q(\hilxmm, p_{m+1}^*(\lmp\tens\dep))=0$$ 
pour $q\ge 2$. On conclut avec le pas (5).$\Box$

\begin{rem}
  {\rm Pour montrer l'annulation de l'espace de cohomologie
  $\H^i(\hilx,\de)$ pour $i>0$, on aurait pu remarquer que le \fib
  canonique de $\hilx$, $\omega_{\hilx}$, est le \fib d{\'e}terminant
  $\deo$ sur $\hilx$ associ{\'e} au
  \fib canonique de $X$, $\ox$, et que l'application de  $\Pic(X)$
  dans $\Pic(\hilx)$ qui associe {\`a} $A$, le \fib d{\'e}terminant associ{\'e} {\`a}
  $A$, $\de$, est un homomorphisme de groupes. Alors l'annulation
  souhait{\'e}e est une cons{\'e}quence du \th de Kawamata-Viehweg, le \fib
  d{\'e}terminant associ{\'e} {\`a} $\ox^{-1}\tens A$ {\'e}tant big et nef. Il n'y a donc pas
  besoin  de toute l'hypoth{\`e}se; l'amplitude  de $\ox^{-1}\tens A$
  suffit pour ce point.}
\end{rem}

{\bf Remerciements~:} Je remercie J. Le Potier pour son enseignement
pr{\'e}cieux, N. Dan pour toute l'aide qu'il m'a apport{\'e}e, C. Mourougane
pour sa pr{\'e}sence, D. Roessler pour la r{\'e}f{\'e}rence \cite{EGA}.



\end{document}